\documentclass[conference]{IEEEtran}
\usepackage{graphicx}
\usepackage{epstopdf}
\usepackage{amsmath,amssymb,epsfig,verbatim}
\usepackage{mathpazo}
\usepackage{microtype}
\usepackage{prettyref}
\usepackage{threeparttable}
\usepackage{multirow}
\usepackage{tabularx}
\usepackage{booktabs}
\usepackage{colortbl}

\usepackage[T1]{fontenc}
\usepackage{textcomp}

\def\F{{\cal F}}

\def\reals{\mathbb{R}}

\def\comp{\raise 1pt \hbox{$\scriptstyle\circ$}}

\def\argmin{\mathop{\rm argmin}\limits}

\def\minimizeQuotes{\mathop{\rm \mbox{``}min\mbox{''}}\limits}
\def\minimize{\mathop{\rm min}\limits}

\def\st{\mathop{\rm subject\ to}}

\def\upto{{\raise 1pt \hbox{$\scriptstyle \,\nearrow\,$}}}
\def\downto{{\raise 1pt \hbox{$\scriptstyle \,\searrow\,$}}}

\def\tos{\rightrightarrows}

\begin{document}
\newcommand{\boldx}{\mbox{$\mathbf{x}$}}
\newcommand{\boldy}{\mbox{$\mathbf{y}$}}
\newcommand{\boldf}{\mbox{$\mathbf{f}$}}
\newcommand{\boldz}{\mbox{$\mathbf{z}$}}
\newcommand{\boldF}{\mbox{$\mathbf{F}$}}
\newcommand{\boldG}{\mbox{$\mathbf{G}$}}
\newcommand{\boldg}{\mbox{$\mathbf{g}$}}
\newcommand{\boldh}{\mbox{$\mathbf{h}$}}
\newcommand{\boldH}{\mbox{$\mathbf{H}$}}
\newcommand{\boldzero}{\mbox{$\mathbf{0}$}}
\newcommand{\Rbb}{\mbox{$\mathbb R$}}

\newcommand{\boldq}{\mbox{$\mathbf{q}$}}
\newcommand{\boldtau}{\mbox{\boldmath$\tau$}}

\newcommand{\pref}{\prettyref}

\newtheorem{definition}{Definition}

\title{Bilevel Optimization based on Iterative Approximation of Multiple Mappings}


\author{Ankur Sinha$^1$, Zhichao Lu$^2$, Kalyanmoy Deb$^2$ and Pekka Malo$^3$
\vspace{2mm}\\
$^1$Production and Quantitative Methods\\
Indian Institute of Management,
Ahmedabad 380015 India\\
{asinha@iima.ac.in}\vspace{2mm}\\
$^2$Department of Electrical and Computer Engineering\\
Michigan State University,
East Lansing, MI, USA\\
{luzhicha@msu.edu, kdeb@egr.msu.edu}\vspace{2mm}\\
$^3$Department of Information and Service Economy\\
Aalto University School of Business,
PO Box 21220, 00076 Aalto, Finland\\
{pekka.malo@aalto.fi}}

\maketitle

\begin{abstract}
A large number of application problems involve two levels of optimization, where one optimization task is nested inside the other. These problems are known as bilevel optimization problems and have been studied by both classical optimization community and evolutionary optimization community. Most of the solution procedures proposed until now are either computationally very expensive or applicable to only small classes of bilevel optimization problems adhering to mathematically simplifying assumptions. In this paper, we propose an evolutionary optimization method that tries to reduce the computational expense by iteratively approximating two important mappings in bilevel optimization; namely, the lower level rational reaction mapping and the lower level optimal value function mapping. The algorithm has been tested on a large number of test problems and comparisons have been performed with other algorithms. The results show the performance gain to be quite significant. To the best knowledge of the authors, a combined theory-based and population-based solution procedure utilizing mappings has not been suggested yet for bilevel problems.
\end{abstract}


\section{Introduction}
\label{intro}
Interest in bilevel optimization has been growing due to a number of new applications that are arising in different fields of science and engineering. Bilevel programming is quite common in the area of defense where these problems are studied as attacker-defender problems. The problem was introduced by Bracken and McGill \cite{bracken73} in the area of mathematical programming, where an inner  optimization problem acts as a constraint to an outer optimization problem. One of the follow-up papers by Bracken and McGill \cite{bracken1974defense} highlighted the applications of bilevel programming in defense. Since then a number of studies on homeland security \cite{brown05,wein09,an13} have been performed, where it is common to have bilevel, trilevel and even multilevel optimization models. In the area of operations research, bilevel optimization is gaining importance in the context of interdiction and protection of hub-and-spoke networks \cite{lei2013identifying}, as most of the critical infrastructures like transportation and communications are predominantly hub-and-spoke. In other game theoretic settings, bilevel optimization has been used in transportation~\cite{migdalas95,constantin95,brotcorne01}, optimal tax policies~\cite{labbe98,my-cec13,my-cec15a}, investigation of strategic behavior in deregulated markets~\cite{hu07}, model production processes~\cite{nicholls95} and optimization of retail channel structures~\cite{williams11}. The applications extend to a variety of other domains, like, facility location \cite{jin2007bi,uno2008evolutionary,sun08}, chemical engineering \cite{smith82,clark1990bilevel}, structural optimization \cite{bendsoe95,christiansen97}, and optimal control \cite{mombaur2010human,albrecht2011imitating} problems. While new applications that are inherently bilevel in nature are arising at a fast pace, the development of computationally efficient algorithms for such problems has not kept the pace with the applications.

A significant body of literature exists on bilevel optimization and its optimality conditions \cite{LiMo02,dempe02,dempe2007new,dempe2014necessary,wiesemann2013pessimistic} in the classical optimization literature. However, on the algorithm front most attention has been given to only simple instances of bilevel optimization where the objective functions and constraints are linear \cite{WeHs91,Be93}, quadratic \cite{BaMo90,EdBa91,AlHoPa92} or convex \cite{liu1998trust}. This is not surprising given the fact that bilevel optimization is difficult to an extent that merely evaluating the bilevel optimality of a given solution is an NP-hard task \cite{vicente94}. Researchers have also attempted to solve these problems using computational techniques like evolutionary algorithms. Most of the bilevel algorithms relying on evolutionary framework have been nested in nature \cite{mathieu,yin-bilevel,li07b,zhu2006hybrid}. One of the drawbacks of such an approach is that it might be able to solve small instances of bilevel problems, but as soon as the problem scales-up beyond a few variables, the computational requirements increase tremendously. However, the evolutionary algorithms still have a niche in solving these problems as it maintains a population at each iteration of the algorithm. A population of points may allow modeling various mappings in bilevel optimization to reduce the computational expense \cite{my-bilevel16}. Some studies in this direction are \cite{my-cec16a,my-ejor17,my-arxiv13,my-cec14}. We believe that exploiting some of the mathematical properties of bilevel problems through modeling of various mappings in bilevel is the way forward in solving such problems.

In this paper, we focus on two important mappings in bilevel optimization borrowed from the mathematical optimization literature. The first mapping is the lower level reaction set mapping (known as $\Psi$-mapping), which provides the lower level optimal solution(s) corresponding to any given upper level vector. Considering the upper level problem as the leader's problem and the lower level problem as the follower's problem, the reaction set mapping represents the rational decisions of the follower corresponding to any decision taken by the leader. The second mapping is the lower level value function mapping (known as $\varphi$-mapping) that provides the optimal objective function value to the follower's problem for any given leader's decision. While the first mapping can be a set-valued mapping, the second mapping is always single-valued. We work with meta-modeling techniques that try to approximate these two mappings and develop a computationally efficient evolutionary algorithm for solving bilevel problems. The algorithm has been tested on a number of test problems, and the computational gain when compared with other techniques is found to be significant. In this paper, we also extend an existing test-suite of bilevel test problems \cite{my-ecj14} with a couple of additional problems to better evaluate our proposed solution procedure.

The paper is organized as follows. To begin with, we provide a brief literature survey of bilevel optimization using evolutionary algorithms. This is followed by various formulations of the bilevel optimization problem and discussion of the two mappings that we approximate in this paper. Thereafter, we provide the bilevel evolutionary optimization algorithm which is an extension of the algorithm proposed in the previous studies \cite{my-ejor17,my-arxiv13,my-cec14}. Following this, we provide the empirical results on a number of test problems. A comparative study with other approaches is also included. Finally, we end the paper with the conclusions section.

\section{A Survey on Evolutionary Bilevel Optimization}
Most of the evolutionary algorithms for bilevel optimization are nested in nature, where one optimization algorithm is used within the other. The outer algorithm handles the upper level problem and the inner algorithm handles the lower level problem. Such a structure necessitates that the inner algorithm is called for every upper level point generated by the outer algorithm. Therefore, nested approaches can be quite computationally demanding, and can only be applied to small scale problems. One can find studies with evolutionary algorithm being used for the upper level problem and classical approach being used for the lower level problem. If the lower level problem is complex, researchers have used evolutionary algorithms at both levels. Below we provide a review of evolutionary bilevel optimization algorithms from the past.


Mathieu et al. \cite{mathieu} was one of the first to propose a bilevel algorithm using evolutionary algorithms. He used a genetic algorithm to handle the upper level problem and linear programming to solve the lower level problem for every upper level member generated using genetic operations. This study was followed by nesting the Frank-Wolfe algorithm (reduced gradient method) within a genetic algorithm in Yin \cite{yin-bilevel}. Other authors utilized similar nested schemes in \cite{li06,li07b,zhu2006hybrid}. Studies involving evolutionary algorithms at both levels include \cite{angelo13,angelo2015study}, where authors have used differential evolution at both levels in the first study, and differential evolution within ant colony optimization in the second study.


Replacing the lower level problem in bilevel optimization with its KKT conditions is a common approach for solving the problem both in classical as well as evolutionary computation literature. Some of the evolutionary studies that utilize this idea include \cite{hejazi2002linear,wang05}. The approach has been popular and even recently researchers are relying on reducing the bilevel problem into single level problem using KKT and solving the reduced problem using evolutionary algorithm, for example, see \cite{wang11,jiang2013application,li2015genetic,wan2013hybrid}.


While KKT conditions can only be applied to problems where the lower level adheres to certain mathematically simplifying assumptions, the researchers are exploring techniques that can solve more general instances of bilevel optimization problems. Some of the approaches are based on meta-modeling the mappings within bilevel optimization, while others may be based on meta-modeling the entire bilevel problem itself. Studies in this direction include \cite{my-ejor17,my-arxiv13,my-cec14}. In this paper, we aim to develop an algorithm that tries to capture two important mappings in bilevel optimization; namely, the lower level reaction set mapping and the lower level value function mapping, in order to reduce the computational complexity of the problem.



\section{Different Bilevel Formulations}\label{sec:sec3}

We will start this section by providing a general formulation for bilevel optimization. This is followed by various proposals that researchers have made for reducing a bilevel problem into a single-level problem. The two levels in a bilevel problem are also known as the leader's (upper) and follower's (lower) problems in the domain of game theory. In general, the variables, objectives and constraints are different for the two levels. The upper level variables are treated as parameters while optimizing the lower level problem.
A general bilevel formulation has been provided below (for brevity, we ignore equality constraints):

\begin{definition}\label{def:bilevel1}
	For the upper-level objective function $F:\reals^n\times\reals^m \to\reals$ and lower-level objective function $f:\reals^n\times\reals^m \to\reals$, the bilevel optimization problem is given by
	\begin{align*}
	\minimizeQuotes_{x_u\in X_U,x_l\in X_L} & F(x_u,x_l)\\
	\st &\\  & \hspace{-12mm}x_l\in \argmin_{x_l \in X_L}
	\lbrace
	f(x_u,x_l) : g_j(x_u,x_l)\leq 0, j=1,\dots,J
	\rbrace\\
	& \hspace{-12mm}G_k(x_u,x_l)\leq 0, k=1,\dots,K
	\end{align*}
	where $G_k:X_U \times X_L \to \reals$, $k=1,\dots,K$ denotes the upper level constraints, and $g_j:X_U\times X_L \to \reals$ denotes the lower level constraints.
\end{definition}


The bilevel problem in Definition \ref{def:bilevel1} is ill-defined if there exists more than one lower level optimal solutions for some upper level variables. In such a case, the decision of the lower level remains unclear. However, in case where the lower level optimal solution is single-valued, the lone optimal solution is the rational choice. Optimizing bilevel problems from either optimistic or pessimistic position are two common approaches that researchers have utilized to handle the ambiguity arising from multiple lower level optimal solutions. In an optimistic position, it is assumed that the lower level chooses that optimal solution which is favorable at the upper level. In a pessimistic position, the upper level optimizes its problem according to the worst case scenario. In other words, the lower level may choose a solution from the optimal set that is least favorable at the upper level. In this paper, we assume an optimistic position while solving bilevel optimization problems.

In case when certain mathematically simplifying assumptions like continuities and convexities are satisfied, often the lower level optimization task in Definition \ref{def:bilevel1} is replaced with its KKT conditions.
However, the reduced formulation is not simple to handle, as it induces non-convexities and discreteness into the problem through the complementary slackness conditions.
We do not utilize any properties of the KKT based reduction in this paper, rather we focus on two different formulations in the development of the evolutionary algorithm in this paper.

\begin{figure*}[htp]
	\begin{minipage}[t]{0.49\linewidth}
		\begin{center}
			\epsfig{file=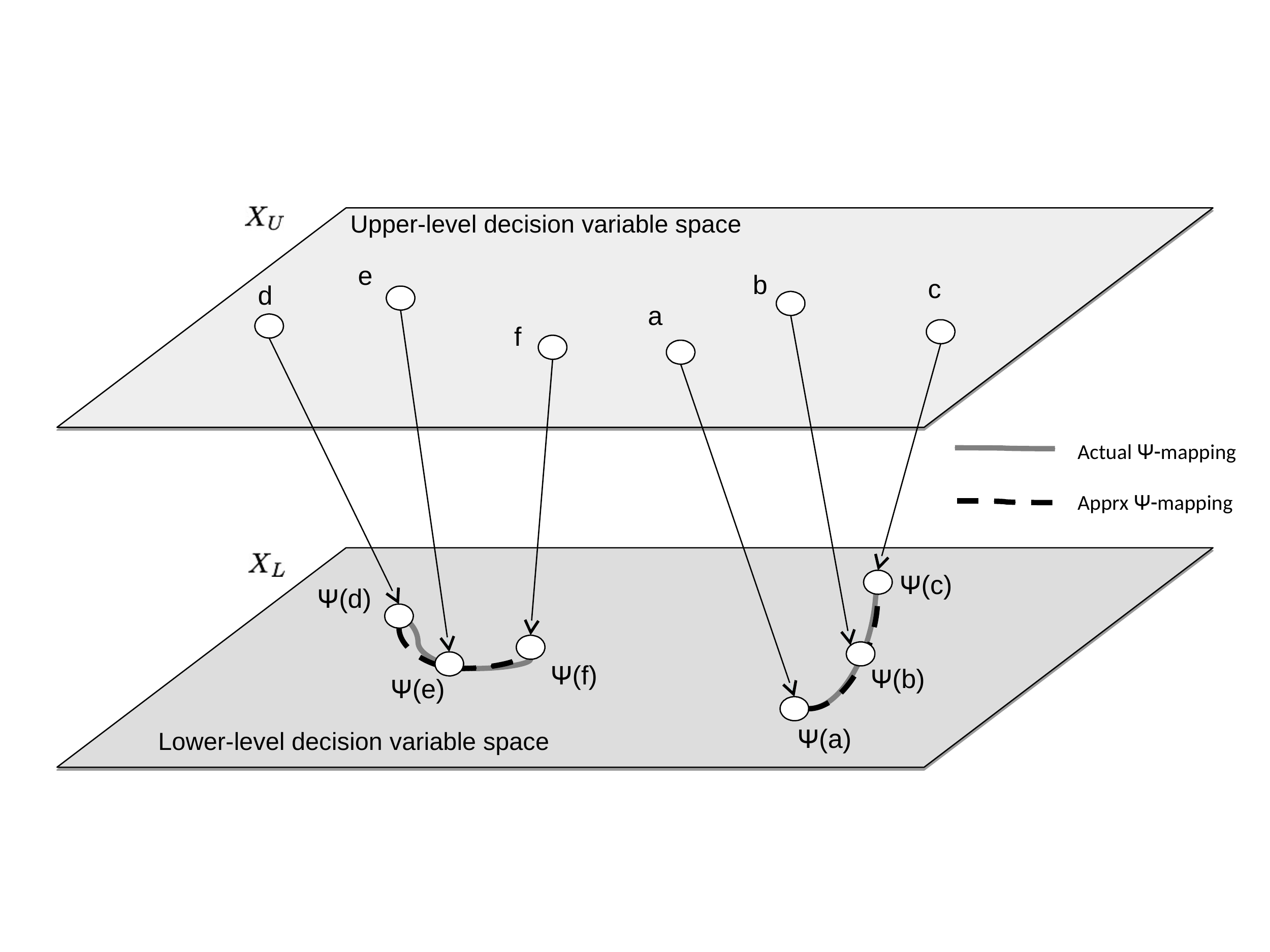,width=\linewidth}
		\end{center}
		\vspace{-12mm}
		\caption{Solving the lower level optimization problem completely for $a, b, c, d, e$ and $f$ provides the corresponding lower level optimal members $\Psi(a), \Psi(b), \Psi(c), \Psi(d), \Psi(e)$ and $\Psi(f)$, where $\Psi$-mapping is assumed to be single valued. Such a mapping can be approximated.}
		\label{fig:explainPsi}
	\end{minipage}\hfill
	\begin{minipage}[t]{0.49\linewidth}
		\begin{center}
			\epsfig{file=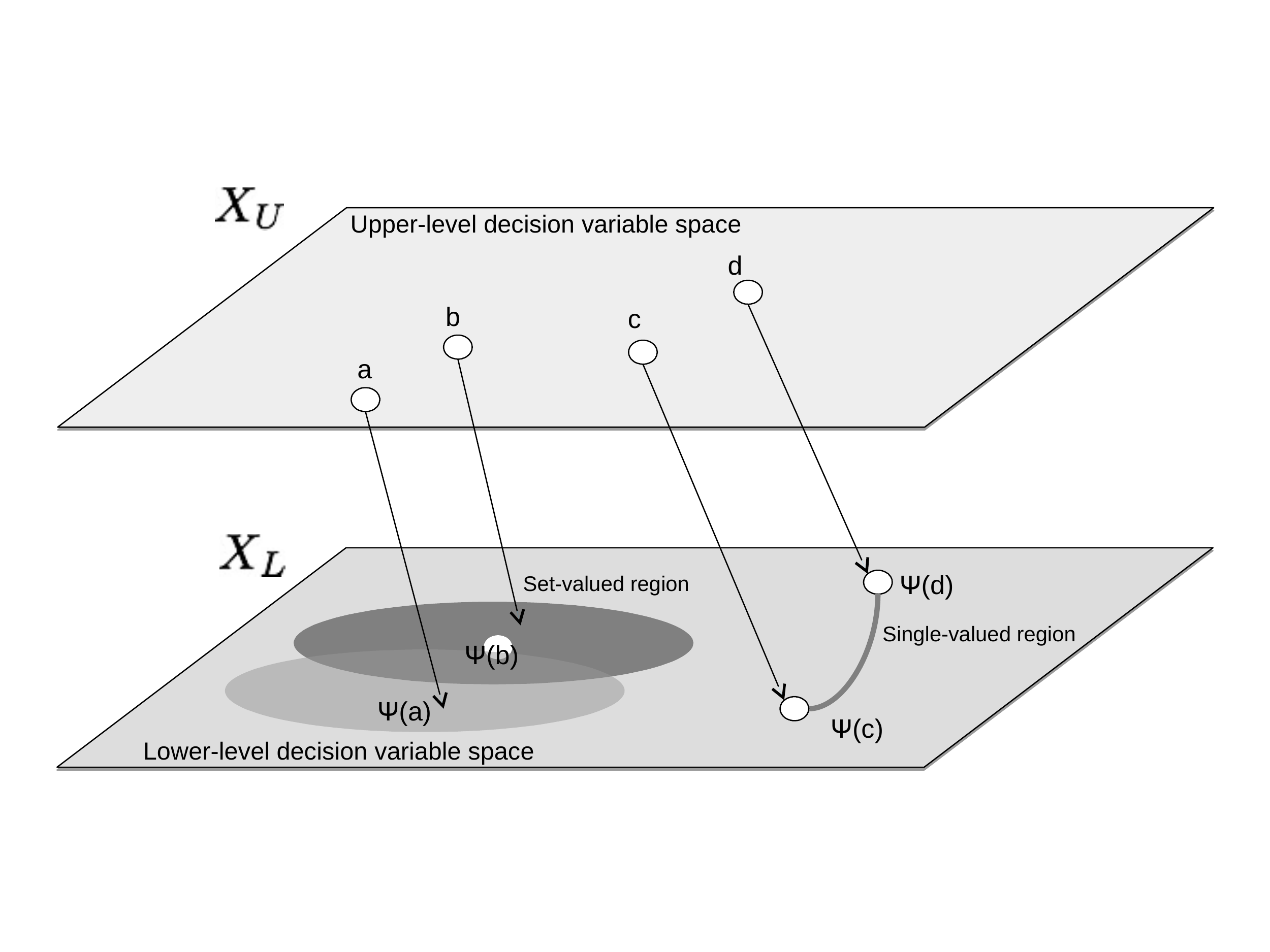,width=\linewidth}
		\end{center}
		\vspace{-12mm}
		\caption{A scenario where the the $\Psi$-mapping is set-valued in some regions and single-valued in other regions. If the $\Psi$-mapping is set-valued then identifying as well as approximating the mapping is not a straightforward task.}
		\label{fig:explainPsiSet}
	\end{minipage}
\end{figure*}

\subsection{Lower Level Reaction Set Mapping}\label{sec:psiMapping}
The formulation provided in Definition \ref{def:bilevel1} can also be stated as follows:
\begin{definition}\label{def:bilevel3}
	Let $\Psi:\reals^n\tos\reals^m$ be the reaction set mapping,
	\begin{align*}
	\Psi(x_u)=&\argmin_{x_l \in X_L}\{f(x_u,x_l) : g_j(x_u,x_l)\leq 0, j=1,\dots,J\},
	\end{align*}
	which represents the constraint defined by the lower-level optimization problem, i.e. $\Psi(x_u)\subset X_L$ for every $x_u\in X_U$. Then the following gives an alternative formulation for the bilevel optimization problem:
	\begin{align*}
	\minimize_{x_u\in X_U, x_l\in X_L}\quad & F(x_u,x_l) \\
	\st\quad  & \\
	& \hspace{-12mm} x_l \in \Psi(x_u) \\
	& \hspace{-12mm} G_k(x_u,x_l)\leq 0, k=1,\dots,K
	\end{align*}
\end{definition}

Using the above definition, a bilevel problem can be reduced to a single level constrained problem given that the $\Psi$-mapping can somehow be determined. Unfortunately this is rarely the case. Studies in the evolutionary computation literature that rely on iteratively approximation of this mapping to reduce the lower level optimization calls could be found in \cite{my-ejor17,my-arxiv13,my-cec14}. To illustrate the idea, let's consider the Figure~\ref{fig:explainPsi}. To acquire sufficient data for constructing the $\Psi$-mapping approximation, a few lower level problems need to be optimized completely for their corresponding upper level decision vectors in the beginning. For instance, the lower level decisions for the upper level decisions $a,b,c,d,e$ and $f$ are determined by optimizing the lower level problem, which are then used to locally approximate the $\Psi$-mapping. This has been shown in Figure~\ref{fig:explainPsi}. Even though the actual $\Psi$-mapping is still unknown, the local approximation can then be substituted to identify the lower level optimal decision for every new upper level member to avoid the lower level optimization task. This procedure of approximating the mapping and utilizing it to predict the lower level optimum needs to be repeated iteratively until convergence to the bilevel optimum. The idea works well when the $\Psi$-mapping is single valued. If the lower level has multiple optimal solutions for some upper level members as shown in Figure \ref{fig:explainPsiSet}, then identifying as well as approximating the mapping is not a straightforward task.

\begin{figure*}[!h]
	\begin{center}
		\begin{minipage}[t]{0.7\linewidth}
			\begin{center}
				\epsfig{file=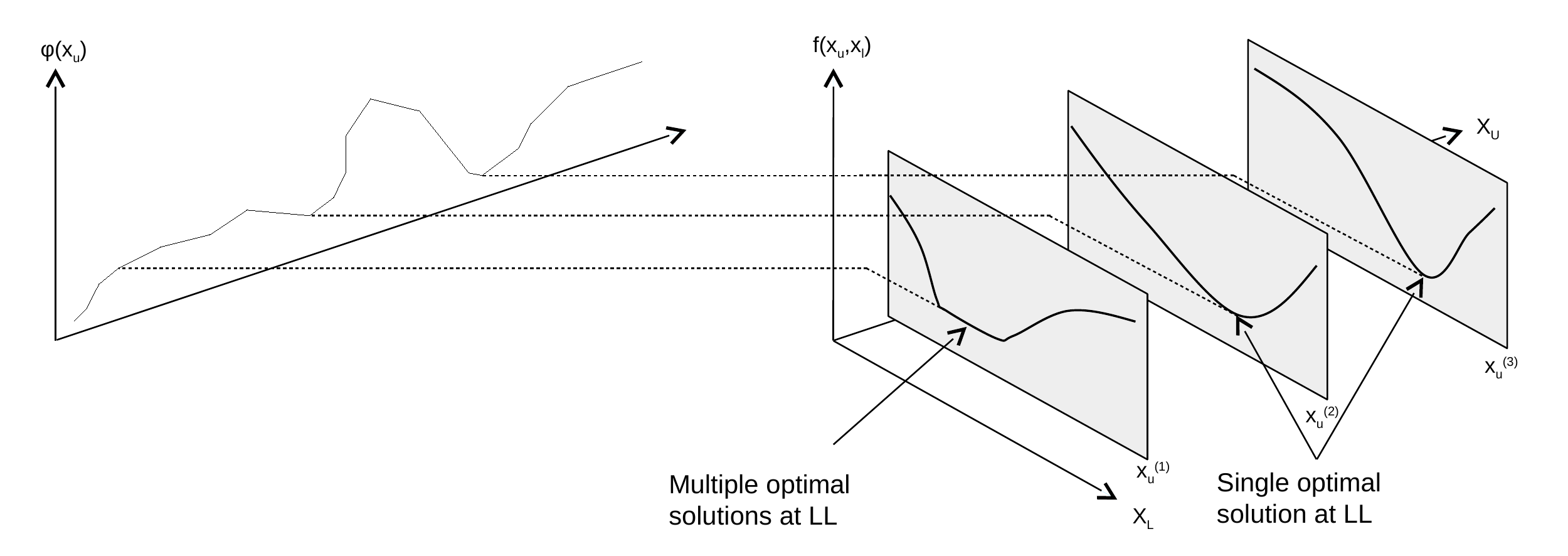,width=\linewidth}
			\end{center}
			\caption{An example showing a $\varphi$-mapping.}
			\label{fig:explain5}
		\end{minipage}
	\end{center}
\end{figure*}

\subsection{Lower Level Optimal Value Function Mapping}\label{sec:varphiMapping}
Another formulation for the bilevel optimization problem in Definition \ref{def:bilevel1} can be written using the optimal lower level value function: \cite{ye2010new}:
\begin{definition}\label{def:bilevel4}
	Let $\varphi: X_U \to R$ be the lower level optimal value function mapping,
	\begin{align*}
	\varphi(x_u)=\minimize_{x_l \in X_L} \{f(x_u,x_l): g_j(x_u,x_l)\leq 0, j=1,\dots,J \},
	\end{align*}
	which represents the optimal function value at the lower level for any given upper level decision vector. Using this lower level optimal value function, the bilevel optimization problem can be expressed as:
	\begin{align*}
	\minimize_{x_u\in X_U, x_l\in X_L}\quad & F(x_u,x_l) \\
	\st\quad  & \\
	& \hspace{-12mm} f(x_u,x_l) \le \varphi(x_u) \\
	& \hspace{-12mm} g_j(x_u,x_l)\leq 0, j=1,\dots,J\\
	& \hspace{-12mm} G_k(x_u,x_l)\leq 0, k=1,\dots,K.
	\end{align*}
\end{definition}

As in the case of $\Psi$-mapping, if the $\varphi$-mapping can be somehow determined, a bilevel problem can be reduced to a single level problem as described in Definition~\ref{def:bilevel4}. Along the process of an algorithm, the $\varphi$-mapping can be approximated and used to solve the reduced single level problem formulation in an iterative manner. Such an evolutionary algorithm has been recently discussed in \cite{my-cec16a}. The approximation of the optimal value function ($\varphi$) mapping is, in general, less complicated than the reaction set ($\Psi$) mapping, in the sense that, the $\varphi$-mapping is always scalar-valued regardless of the lower level variable dimension and whether or not there exist multiple lower level optimal solutions. However, the $\varphi$-mapping based reduction is not necessarily always better than the $\Psi$-mapping based reduction. Definition \ref{def:bilevel4} requires the problem to be solved with respect to upper as well as lower level variables, while in Definition \ref{def:bilevel3} the lower level variables are readily available from the $\Psi$-mapping. The $\Psi$-mapping based reduction also contains fewer constraints. Therefore, clearly there is a trade-off.

\section{Evaluating the Performance of $\Psi$ and $\varphi$ mappings on test problems}
\begin{figure*}[htp]
	\begin{center}
		\epsfig{file=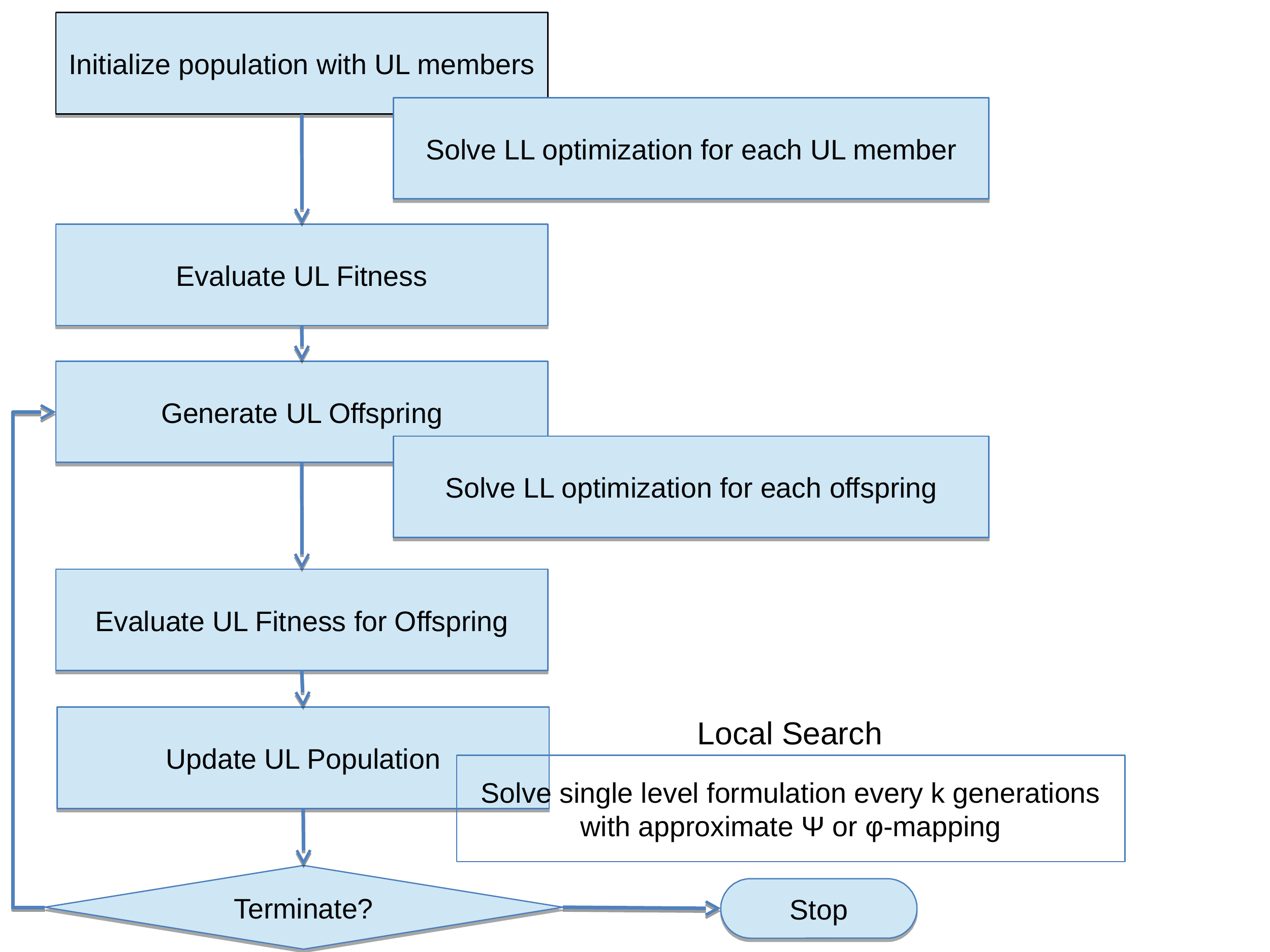,width=0.60\linewidth}
		\caption{Nested approach with evolutionary algorithm at upper level (UL) and SQP at lower level (LL). Local search based on $\Psi$ or $\varphi$ mapping may be performed to make the nested approach faster.}
		\label{fig:nested}
	\end{center}
\end{figure*}

In this section, we implement the $\Psi$ and $\varphi$ mappings separately in two different nested algorithms to evaluate the advantages and disadvantages of using the two mappings as a local search. For evaluating the two mappings, we choose a set of simple test problems that are provided in Tables~\ref{tab:tpTable1} and~\ref{tab:tpTable2}. Firstly, we create a nested algorithm that utilizes an evolutionary approach for solving the upper level problem and sequential quadratic programming (SQP) for solving the lower level problem.  Most of the lower level problems in the considered test cases being convex, explains the choice for using sequential quadratic programming (SQP) at the lower level. We enhance the nested approach by allowing it to approximate the $\Psi$ and $\varphi$ mappings and measure the performance gain provided by using each of the mappings separately. The implementation of the approaches has been outlined through the Figure \ref{fig:nested}. The flowchart without the overlapping box provides the steps involved in the nested approach. In case the idea involving $\Psi$ and $\varphi$ mappings has to be used then the local search (as mentioned in the overlapping box) is conducted every $k$ generations of the nested algorithm after the update step. A detailed description of the nested algorithm has been provided below.
\begin{enumerate}
\item Create a random population of size N comprising of upper level variables
\item Solve the lower level optimization problem using SQP for each upper level variable.
\item Evaluate the fitness of each population member using upper level function and constraints (refer to Section \ref{sec:fitnessAssignment})
\item Choose $2\mu$ population members using tournament selection and apply genetic operators (refer to Section \ref{sec:geneticOperations}) to produce $\lambda$ offspring.
\item Solve the lower level optimization problem using SQP for each offspring.
\item Evaluate the fitness of each offspring using using upper level function and constraints
\item Form a pool consisting of $r+\lambda$ members, where $r$ members are chosen randomly from the population and $\lambda$ members are the offspring. Use the best $r$ members from this pool to replace the chosen $r$ members from the population.
\item Perform a termination check (refer to Section \ref{sec:termination}) and proceed to Step 5 if termination check is false, otherwise stop.
\end{enumerate}
The parameters used in the implementation of the above procedure are $N=50, \mu=2, \lambda=3$ and $r=2$. The lower level SQP terminates when the improvement in the lower level function value is less than 1e-6.

\subsection{Approximating the $\Psi$-mapping}\label{sec:psiApprox}
Let $\mathcal{H}$ be the hypothesis space. The hypothesis space consists of all functions that can be used to generate a mapping between the upper level decision vectors and optimal lower level decision vectors. Given a sample $\mathcal{S}$, consisting of upper level points and corresponding optimal lower level points, we would like to identify a model $\hat{\Psi}\in\mathcal{H}$ that minimizes the empirical error on the sample, i.e.
\begin{equation}\label{eq:empirical}
\hat{\psi}=\argmin_{h \in \mathcal{H}}\sum_{i\in \mathcal{I}} L(h(x_{u}^{(i)}),\bar{x}_{l}^{(i)}),
\end{equation}
where $L:X_L\times X_L\to\reals$ denotes the prediction error, $x_{u}^{(i)}$ is any given upper level vector and $\bar{x}_{l}^{(i)}$ is its corresponding optimal solution. The prediction error may be calculated as follows:
$$
L(h(x_{u}^{(i)}),\bar{x}_{l}^{(i)})=|\bar{x}_{l}^{(i)}-h(x_{u}^{(i)})|^2.
$$
We have restricted the hypothesis space $\mathcal{H}$ to consist of second-order polynomials which reduces the error minimization problem to an ordinary quadratic regression problem. The sample $\mathcal{S}$ can be created from the population members or an archive. It should be noted that this can approximate only single-valued mapping and will fail if the mapping becomes set-valued.

\subsection{Approximating the $\varphi$-mapping}\label{sec:varphiApprox}
Once again, let $\mathcal{H}$ be the hypothesis space of functions and $\mathcal{S}$ be a sample of upper and corresponding lower level points, our aim is to identify a model $\hat{\varphi}\in\mathcal{H}$ that minimizes the empirical error on the sample, i.e.
\begin{equation}\label{eq:empirical}
\hat{\varphi}=\argmin_{u \in \mathcal{H}}\sum_{i\in \mathcal{I}} L(u(x_{u}^{(i)}),\bar{f}^{(i)}),
\end{equation}
where $L:\reals \times \reals \to\reals$ denotes the prediction error, $x_{u}^{(i)}$ is any given upper level vector and $\bar{f}^{(i)}$ is its corresponding optimal function value. The prediction error can once again be computed as follows:

$$
L(u(x_{u}^{(i)}),f^{(i)})=|\bar{f}^{(i)}-u(x_{u}^{(i)})|^2.
$$
We have once again restricted the hypothesis space $\mathcal{H}$ to consist of second-order polynomials. Since the $\varphi$-mapping is always single valued, approximating it will not involve similar issues as for the $\Psi$-mapping.

\section{Comparison Results for $\Psi$ vs $\varphi$ Approximations}
\begin{figure*}
	\begin{minipage}[t]{0.49\linewidth}
		\begin{center}
			\epsfig{file=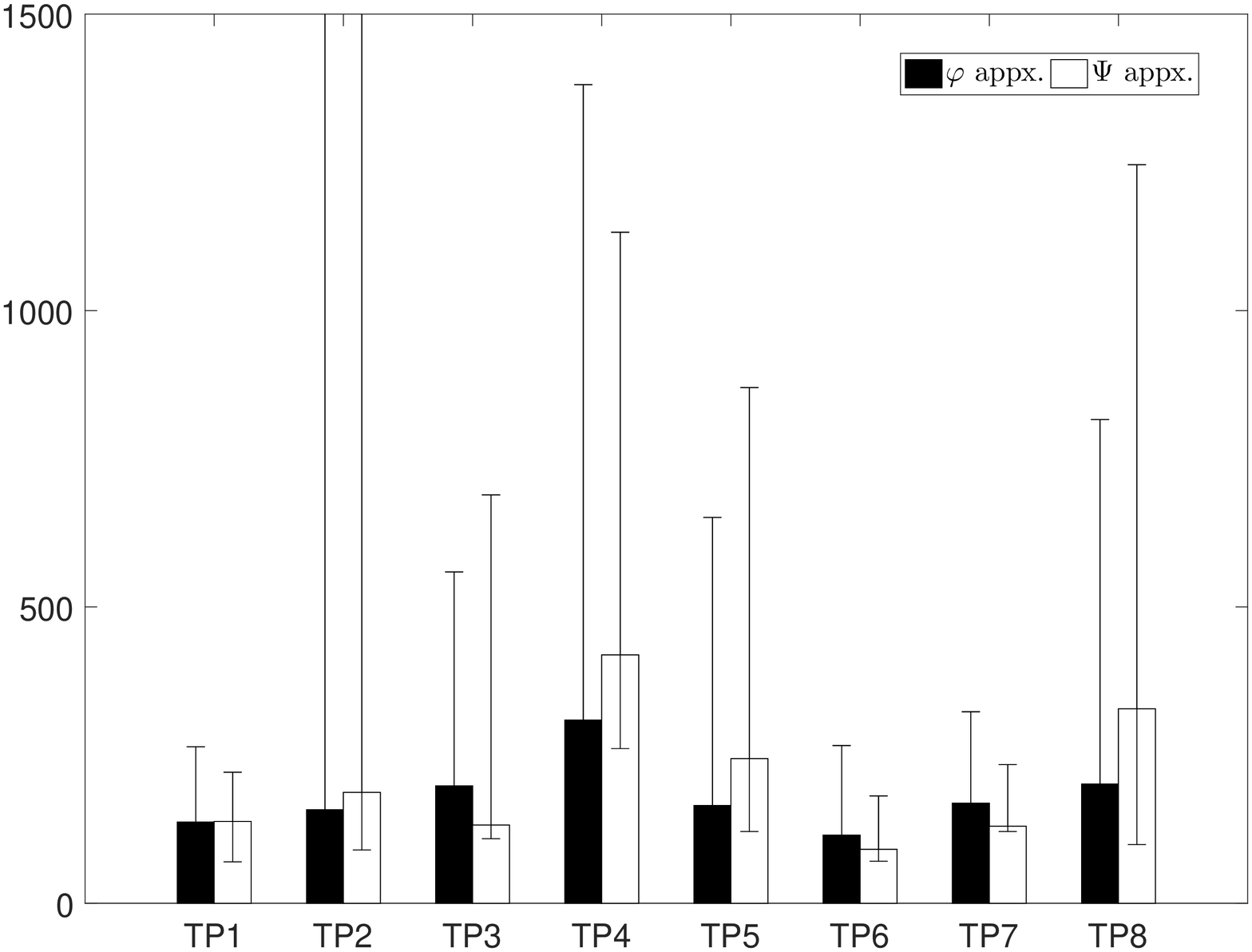,width=0.9\linewidth}
		\end{center}
		\caption{Error plot from 31 runs for the upper level function evaluations on test problems 1 to 8.}
		\label{fig:resFigure1}
	\end{minipage}\hfill
	\begin{minipage}[t]{0.49\linewidth}
		\begin{center}
			\epsfig{file=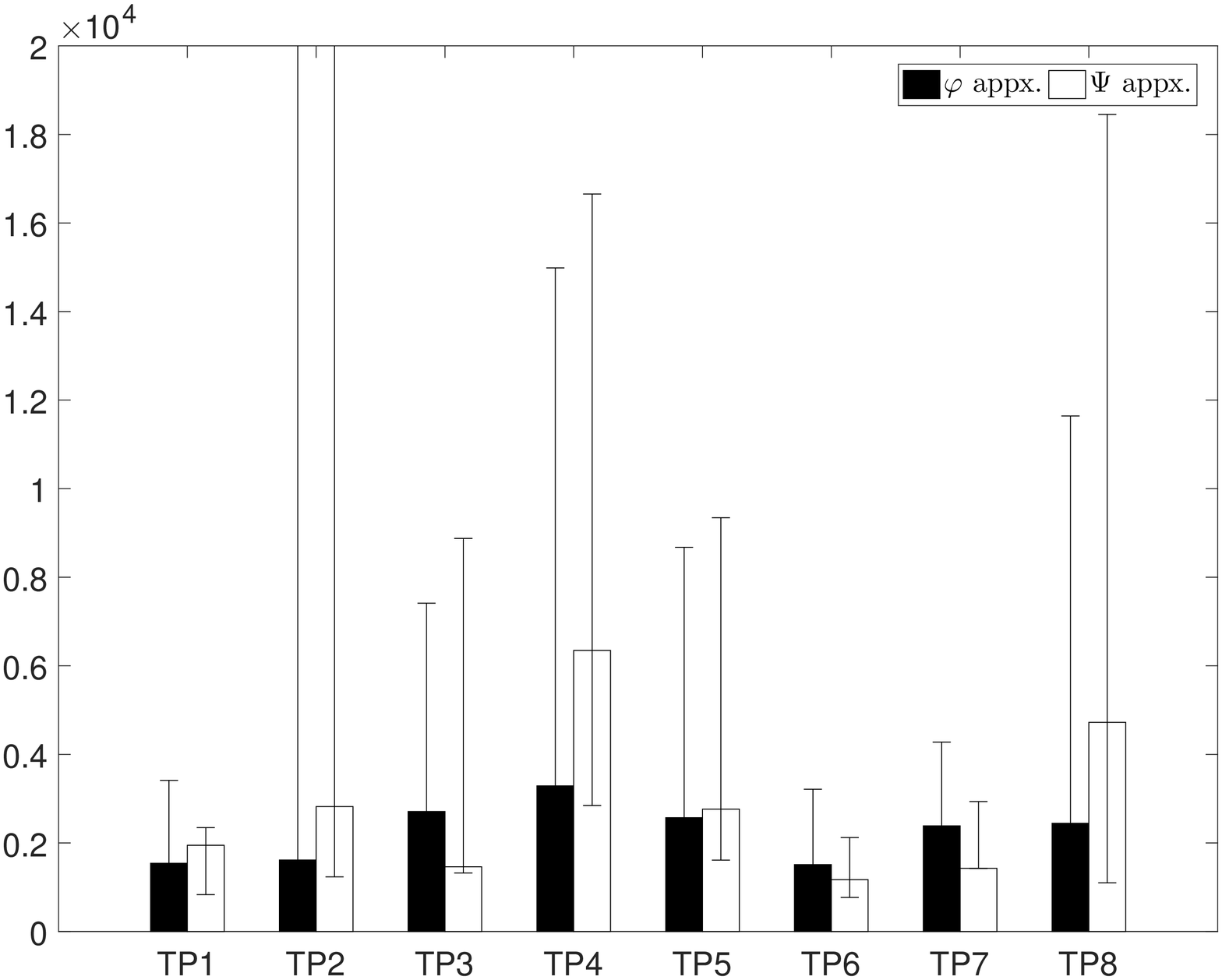,width=0.9\linewidth}
		\end{center}
		\caption{Error plot from 31 runs for the lower level function evaluations on test problems 1 to 8.}
		\label{fig:resFigure2}
	\end{minipage}
\end{figure*}

For comparing $\Psi$-approximation approach against $\varphi$-approximation approach, we use a set of 8 test problems selected from the literature given in Tables \ref{tab:tpTable1} and \ref{tab:tpTable2}. Table \ref{tab:resTable1} compares the median function evaluations at both level for three algorithms $\Psi$-approximation, $\varphi$-approximation and nested algorithm. The results have been produced from 31 runs of the algorithm and further details about the runs can be found in Figures \ref{fig:resFigure1} and \ref{fig:resFigure2}. Both $\Psi$-approximation and  $\varphi$-approximation perform equally well and outperform the nested approach in this study. The differences in the performance of $\Psi$-approximation and  $\varphi$-approximation can be attributed to differences in the quality of approximation produced during the intermediate steps of the algorithm.
In Table~\ref{tab:resTable1-2}, we provide a comparison of the meta-modeling results with other evolutionary approaches \cite{wang05,wang11} to provide an idea about the extent of savings that can be produced using meta-modeling techniques. The advantage is quite clear as the savings are better by multiple order of magnitudes on the set of test problems considered in this study.


\begin{table*}[]
	\centering
	\caption{Median function evaluations required at upper level (UL) and the lower level (LL) from 31
		runs of $\Psi$-approximation algorithm, $\varphi$-approximation algorithm and nested algorithm.}
	\vspace{0mm}
	\label{tab:resTable1}
	\begin{tabular}{@{}ccccccccc@{}}
		\toprule
		& \multicolumn{3}{c}{UL Func. Evals.} & \multicolumn{3}{c}{LL Func. Evals.}      & \multicolumn{2}{c}{Savings}   \\ \midrule
		& $\varphi$ Appx & $\Psi$ Appx & Nested & $\varphi$ Appx & $\Psi$ Appx & Nested & $\frac{\mbox{Nested}-\mbox{$\varphi$ Appx}}{\mbox{Nested}}$ & $\frac{\mbox{Nested}-\mbox{$\Psi$ Appx}}{\mbox{Nested}}$ \\
		& Med            & Med         & Med     & Med            & Med         & Med     &           &        \\ \midrule
TP1  	&	137	&	138	&	-	&	1539	&	1948	&	-	&	Large	&	Large	\\	
TP2  	&	158	&	187	&	444	&	1614	&	2820	&	5252	&	69\%	&	47\%	\\	
TP3  	&	198	&	132	&	642	&	2710	&	1461	&	6530	&	59\%	&	78\%	\\	
TP4  	&	309	&	419	&	1760	&	2976	&	6347	&	18073	&	83\%	&	66\%	\\	
TP5  	&	165	&	244	&	633	&	2571	&	2763	&	6616	&	62\%	&	59\%	\\	
TP6  	&	115	&	91	&	153	&	1512	&	1170	&	2169	&	30\%	&	46\%	\\	
TP7  	&	169	&	130	&	196	&	2386	&	1426	&	2606	&	9\%	&	44\%	\\	
TP8  	&	201	&	328	&	423	&	2443	&	4722	&	8233	&	69\%	&	42\%	\\	\bottomrule
	\end{tabular}
\end{table*}

\begin{table}[]
	\centering
	\caption{Mean of total function evaluations (UL evaluations +LL evaluations) required by different approaches.}
	\vspace{0mm}
	\label{tab:resTable1-2}
	\begin{tabular}{@{}cccccc@{}}
		\toprule
		& \multicolumn{5}{c}{Mean Func. Evals. (UL+LL)}     \\ \midrule
		& $\varphi$ Appx. & $\Psi$ Appx. & Nested & WJL \cite{wang05}     & WLD \cite{wang11}    \\ \midrule
		TP1 & 1,611      & 2,421   & 34,462     & 85,499   & 86,067  \\
		TP2 & 1,923      & 3,262   & 6,235      & 256,227  & 171,346 \\
		TP3 & 2,624      & 1,482   & 8,125      & 92,526   & 95,851  \\
		TP4 & 3,612      & 6,721   & 19,948     & 291,817  & 211,937 \\
		TP5 & 2,812      & 3,388   & 7,398      & 77,302   & 69,471  \\
		TP6 & 1,578      & 1,034   & 1,534      & 163,701  & 65,942  \\
		TP7 & 2,110      & 1,456   & 2,286      & 1,074,742 & 944,105 \\
		TP8 & 2,734      & 4,434   & 5,325      & 213,522  & 182,121 \\ \bottomrule
	\end{tabular}
\end{table}

\begin{table}[]
	\centering
	\caption{Statistics for upper level function evaluations for $\varphi$-approximation algorithm on the modified test problems (m-TP). The $\varphi$-approximation algorithm and nested algorithm fail on all the modified test problems.}
	\vspace{0mm}
	\label{tab:resTable2}
	\begin{tabular}{@{}cccccc@{}}
		\toprule
		& \multicolumn{3}{c}{$\varphi$ Appx.}      & $\Psi$ Appx. & Nested \\ \midrule
		& Min         & Med  & Max   & Min/Med/Max & Min/Med/Max \\ \midrule
		m-TP1  & 138         & 192 & 344  & -           & -          \\
		m-TP2  & 124         & 236 & -    & -           & -          \\
		m-TP3  & 140         & 242 & 699  & -           & -          \\
		m-TP4  & 185         & 545 & 2,582 & -           & -          \\
		m-TP5  & 172         & 242 & 977  & -           & -          \\
		m-TP6  & 159         & 181 & 559  & -           & -          \\
		m-TP7  & 119         & 227 & 501  & -           & -          \\
		m-TP8  & 158         & 462 & 2,119 & -           & -          \\ \bottomrule
	\end{tabular}
\end{table}

\begin{table}[]
	\centering
	\caption{Statistics for lower level function evaluations from 31
		runs of the $\varphi$-approximation algorithm on the modified test problems (m-TP). The $\varphi$-approximation algorithm and nested algorithm fail on all the modified test problems.}
	\vspace{0mm}
	\label{tab:resTable3}
	\begin{tabular}{@{}cccccc@{}}
		\toprule
		& \multicolumn{3}{c}{$\varphi$ Appx.}      & $\Psi$ Appx. & Nested \\ \midrule
		& Min         & Med  & Max   & Min/Med/Max & Min/Med/Max \\ \midrule
		m-TP1  & 1,988        & 2,477 & 8,334  & -           & -           \\
		m-TP2  & 2,394        & 4,420 & -     & -           & -           \\
		m-TP3  & 1,404        & 3,321 & 12,353 & -           & -           \\
		m-TP4  & 1,911        & 5,632 & 25,356 & -           & -           \\
		m-TP5  & 3,129        & 4,166 & 15,345 & -           & -           \\
		m-TP6  & 2,498        & 3,464 & 9,325  & -           & -           \\
		m-TP7  & 1,476        & 5,635 & 12,256 & -           & -           \\
		m-TP8  & 2,721        & 6,324 & 28,993 & -           & -           \\ \bottomrule
	\end{tabular}
\end{table}

It should be noted that the $\Psi$ approximation idea would fail if the $\Psi$-mapping in bilevel optimization is set valued. Next, we test this hypothesis, by modifying the 8 test problems such that each test problem necessarily has a set-valued $\Psi$-mapping. To achieve this, we add two additional lower level variables ($y_p$ and $y_q$) in each test problem. Both the upper and lower level functions are modified as shown below:
\begin{align*}
F^{new}(x_{u},x_{l}) & = F(x_{u},x_{l})+y_{p}^{2}+y_{q}^{2}\\
f^{new}(x_{u},x_{l}) & = f(x_{u},x_{l})+(y_{p} - y_{q})^{2}\\
y_{p}, y_{q} & \in [-1,1]
\end{align*}
The modification makes the lower level problem have infinitely many optimal solutions (for all $y_{p} = y_{q}$) for any given upper level vector. Out of the many optimal solutions the upper level prefers the solution where $y_{p} = y_{q} = 0$. After this simple modification, we once again solve the test problems using $\varphi$-approximation and $\Psi$-approximation approaches. As shown in Tables \ref{tab:resTable2} and \ref{tab:resTable3}, the $\varphi$-approximation algorithm still works but $\Psi$-approximation algorithm completely fails. Function evaluations for $\varphi$-approximation algorithm increases slightly than before because of additional variables in the problem.


Therefore, the $\varphi$-approximation idea clearly has an advantage over the $\Psi$-approximation idea. Moreover, $\varphi$-mapping is always a scalar valued mapping as compared to $\Psi$-mapping which is usually vector valued and can also be set-valued.
However, there is a trade-off. The reduced single level problem formed using $\Psi$-mapping may usually be a little easier to handle as compared to the single level problem formed using $\varphi$-mapping. The reason being that in case of $\Psi$-mapping the lower level variables are readily available, and the reduced problem does not involve lower level constraints. For $\varphi$-mapping, the reduced problem has to be solved both with respect to upper and lower level variables, and the formulation involves both upper and lower level constraints. Given the pros and cons of using the two mappings, next, we would like to develop an evolutionary algorithm that is capable of utilizing the better of the two mappings while solving a bilevel optimization problem.


\section{Bilevel Evolutionary Algorithm based on $\Psi$ and $\varphi$-mapping Approximations}
In this section, we provide the bilevel evolutionary algorithm that approximates the $\Psi$ as well as the $\varphi$ mapping during the intermediate steps of the algorithm. From the previous experiments and the properties of the two mappings we infer that there can be situations when the approximation of the $\Psi$-mapping may fail, while when $\Psi$-mapping can be approximated it offers the advantage of completely ignoring the lower level functions and constraints. Acknowledging this fact, we utilize both the approximations in our algorithm. The algorithm adaptively decides to use one of the mappings based on the quality of fit obtained when approximating the two mappings. Local quadratic approximations are created for the two mappings from a sample of points in the vicinity of the point around which we want to create an approximation. Introducing local approximation, is expected to improve the quality of approximations significantly. Given a sample dataset, the steps for creating an approximation are the same as discussed in Sections~\ref{sec:psiApprox} and~\ref{sec:varphiApprox}. The algorithm also maintains an archive so as to maintain a large dataset for creating and validating the approximations. Deviating from the nested algorithm, we employ the approximated $\Psi$ and $\varphi$ mappings to avoid frequent lower level optimization calls. An earlier version of the algorithm \cite{my-ejor17,my-arxiv13,my-cec14} that relied on $\Psi$-mapping approximation alone was referred as Bilevel Evolutionary Algorithm based on Quadratic Approximations (BLEAQ). We keep the same terminology and refer to the newer version of the algorithm as BLEAQ-II. The pseudocode for the algorithm has been provided in Table \ref{tab:pseudocode}.

\begin{table*}[htp]
	{\small
		\caption{Step-by-step procedure for BLEAQ-II}\label{tab:BLEAQ-II}
		\begin{center}
			\begin{tabular}{p{.03\textwidth}  p{.92\textwidth}}
				Step & Description \\
				\toprule
				{\bf 1} & {\bf Initialization}: \\
				& Generate an initial upper level population $x_u^{(1)},\dots,x_u^{(N)}$ randomly or by a problem-specific method (see Section \ref{sec:initialization}).\\ \\
				& {\bf (a)} For each $x_u^{(j)}$, find a corresponding optimal lower level solution $x_l^{(j)}\in\Psi(x_u^{(j)})$ by solving the lower level problem. Set $\mathcal{P}=\{(x_u^{(j)},x_l^{(j)}),j=1,\dots,N\}$  (see Section~\ref{sec:lowerLevelAlgorithm}).\\ \\
				& {\bf (b)} Tag all vectors $(x_u^{(j)},x_l^{(j)})\in\mathcal{P}$ for which a lower level optimization has been successfully performed as {\bf 1} and store them in the archive $\mathcal{A}$. \\ \\
				& {\bf (c)} Assign fitness to all the members based on upper level function and constraints (refer to Section \ref{sec:fitnessAssignment}). \\
				
				\midrule
				{\bf 2} & {\bf Reproduction}:\\
				& {\bf (a)} {\bf Parent selection}: Randomly choose $2 \mu$ members from the population $\mathcal{P}$, and perform a tournament selection based on the upper level fitness. This produces $\mu$ parents, denoted by $\mathcal{P}_{\text{par}}$.\\ \\
				& {\bf (b)} {\bf Offspring generation}: Create $\lambda$ offsprings, denoted by $\mathcal{P}_{\text{off}}$, from the set of parents $\mathcal{P}_{\text{par}}$  using genetic operators (refer to Section \ref{sec:geneticOperations}).\\
				
				\midrule
				{\bf 3} & {\bf Offspring Update}: \\
				& For each offspring $x^{(j)}=(x^{(j)}_u,x^{(j)}_l)\in\mathcal{P}_{\text{off}}$ produced in the previous step, update the lower level decision $x^{(j)}_l$  using the following strategy:\\ \\
				& {\bf (a)} {\bf Optimization}: If the number of Tag {\bf 1} members in $\mathcal{P}$ is less than half of the size of $\mathcal{P}$, perform lower level optimization to ensure that $x_l^{(j)}\in{\Psi}(x_u^{(j)})$ (as described in Step 1.(b)). If the lower level optimization is successful, tag the offsprings as {\bf 1} and add it to the archive $\mathcal{A}$. \\ \\
				& {\bf (b)} {\bf Approximations}: If the number of Tag {\bf 1} members in $\mathcal{P}$ is more than half of the size of $\mathcal{P}$, for each offspring $x^{(j)}=(x^{(j)}_u,x^{(j)}_l)\in\mathcal{P}_{\text{off}}$, use its neighboring members in the archive $\mathcal{A}$ to construct a local quadratic approximation for $\Psi$-mapping ($q_{\Psi}$) as well as $\varphi$-mapping ($q_{\varphi}$). Compare the mean squared error of the approximations ($e_{mse}^{\Psi}$, $e_{mse}^{\varphi}$). If $e_{mse}^{\Psi} \leq e_{mse}^{\varphi}$ then update the lower level decision associated with the upper level $x_u^{(j)}$ by setting $x^{(j)}_l=q_{\Psi}(x_u^{(j)})$; otherwise solve the auxiliary optimization problem in Section \ref{sec:intermediatePhi} by fixing $x^{(j)}_u$ and varying $x^{(j)}_l$; the optimal $x^{(j)}_l$ is paired with $x^{(j)}_u$ to form the offspring.\\
				
				\midrule
4 & {\bf Improvements}: \\
				& Identify the Tag {\bf 1} member in the current generation in $\mathcal{P}$ with the best fitness, denoted as $x_{best}^{(j)}$. Perform a local search in the vicinity of $x_{best}^{(j)}$ after every $k$ generations and update $x_{best}^{(j)}$ if there is an improvement.   \\ \\
				& {\bf (a)} {\bf Local search}: Construct local quadratic approximations of both $\Psi$-mapping and $\varphi$-mapping using members in the vicinity of $x_{best}^{(j)}$ in the archive $\mathcal{A}$ and record the mean squared error of the approximations ($e_{mse}^{\Psi}$, $e_{mse}^{\varphi}$). Apply local search in $x_{best}^{(j)}$ vicinity using one of the two single level reduction methods described in Sections \ref{sec:psiMapping} and \ref{sec:varphiMapping} (refer to Section \ref{sec:localSearch}).\\ \\
				
				\midrule
				5 & {\bf Termination check:} \\
				& Perform a termination check. If false, proceed to the next generation (Step 2). \\
				\bottomrule			
			\end{tabular}
		\end{center}
		\label{tab:pseudocode}
	}
\end{table*}%

\begin{table*}[htp]
	{\small
		\caption{The lower level evolutionary algorithm is described below that takes an upper level member as input and solves the corresponding lower level problem.}\label{tab:lowerLevelAlgorithm}
		\begin{center}
			\begin{tabular}{p{.03\textwidth}  p{.92\textwidth}}
				Step & Description \\
				\toprule
				{\bf 1} & {\bf Initialization}: Generate an initial lower level population $x_l^{(1)},\dots,x_l^{(N)}$ randomly and assign fitness using lower level objective and constraints.\\ \\
				{\bf 2} & {\bf Genetic Operations}: Randomly choose $2\mu$ members from the population, and perform a tournament selection leading to $\mu$ parents. Create $\lambda$ offspring using the genetic operations described in Section \ref{sec:geneticOperations}. Assign fitness to each offspring. \\ \\
				{\bf 3} & {\bf Update}: Choose $r$ members randomly from the population and pool them with $\lambda$ offspring. Sort the pool by fitness and replace the $r$ members from the population by the best $r$ members from the pool.\\ \\
				{\bf 4} & {\bf Termination check:} Perform a termination check as described in Section \ref{sec:termination}. If false, proceed to the next generation (Step 2). \\ \\
				\bottomrule			
			\end{tabular}
		\end{center}
		\label{pseudocode}
	}
\end{table*}%

\subsection{Initialization}\label{sec:initialization}
The initialization in the algorithm is done by creating random upper level members $x_u^{(1)},\dots,x_u^{(N)}$, and then solving the lower level optimization problem for each member to get optimal $x_l^{(1)},\dots,x_l^{(N)}$. There can be situations, where finding random feasible upper and lower level pair that satisfy both lower and upper level constraints in the problem can be difficult. In such situations, one can solve the following problem to create $(x_u^{(i)},x_l^{(i)})$ pairs that satisfies all the constraints to begin with.
\begin{align*}
\minimize_{x_u\in X_U, x_l\in X_L}\quad & 0 \\
\st\quad  & \\
 & \hspace{-12mm} G_k(x_u,x_l)\leq 0, k=1,\dots,K,\\
 & \hspace{-12mm} g_j(x_u,x_l)\leq 0, j=1,\dots,J.\\
\end{align*}
The above problem can be solved using any standard procedure like a greedy GA or SQP with a random starting point to arrive at a feasible solution. As soon as a feasible member is found, the method stops. Solving the above method repeatedly with random population (in case of GA) or a random starting point (in case of SQP) can provide the starting population of upper level members $x_u^{(1)},\dots,x_u^{(N)}$ for the BLEAQ-II algorithm. For this given set of upper level members, we know that at least one feasible lower level member exists and we still need to solve the lower level problem to find the optimal lower level solutions $x_u^{(1)},\dots,x_u^{(N)}$. \footnote{In case of upper level constraints containing both upper and lower level variables, one can find it difficult to arrive at a $(x_u^{(i)},x_l^{(i)})$ pair that is feasible with respect to all the constraints and the lower level vector is optimal for the given upper level vector. Many formulations of bilevel optimization, therefore, do not consider lower level variables in upper level constraints.}

\subsection{Constraint handling and Fitness Assignment}\label{sec:fitnessAssignment}
The proposed approach always assigns higher fitness to a feasible member over a non-feasible member. For two given members, $(x_u^{(i)},x_l^{(i)})$ and $(x_u^{(j)},x_l^{(j)})$, if both members are feasible with respect to constraints then it looks at the function value. If both members are infeasible, then it looks at the overall constraint violation. This fitness assignment scheme is similar to the one proposed in \cite{debpenalty}. At the lower level the idea can be implemented directly using lower level constraints and lower level function value. At the upper level, for any given upper and lower level pair, we only consider upper level function and constraints, without considering if the corresponding lower level vector is optimal. The information about a lower level vector corresponding to an upper level vector being optimal is stored using tagging (0 or 1).

\subsection{Genetic operations}\label{sec:geneticOperations}
Offspring are produced in the BLEAQ-II approach using standard crossover and mutation operators. Genetic operations at the upper level involve only upper level variables, and the operations at the lower level involve only lower level variables. We utilize parent centric crossover (PCX) and polynomial mutation for generating the offspring. The crossover operator used in the algorithm is similar to the parent-centric crossover (PCX) operator proposed in \cite{my-cec06}. The operator uses three parents and produces offspring around the index parent as described below.
\begin{equation}
c = z^{(p)} + \omega_{\xi}d + \omega_{\eta}\frac{p^{(2)}-p^{(1)}}{2}
\label{eq:child}
\end{equation}
where,
\begin{itemize}
	\item $z^{(p)}$ is the {\em index\/} parent (the best parent among three parents)
	\item $d=z^{(p)}-g$, where $g$ is the mean of $\mu$ parents
	\item $p^{(1)}$ and $p^{(2)}$ are the other two parents
	\item $\omega_{\xi}=0.1$ and $\omega_{\eta}=0.1$ are the two parameters.
\end{itemize}

\subsection{Termination Criteria}\label{sec:termination}
A variance based termination criterion has been used at both levels; some other termination criterion like termination based on no improvement may also be used. Variance based termination allows the algorithm to terminate automatically when the variance of the population becomes small. At the upper level, the variance of the population at any generation, $T$, is computed as follows:
\begin{equation}
\begin{array}{l}
	\alpha_u^{T} = \frac{\sum_{i=1}^{n} \sigma^2(x_{i})|_{T}}{\sum_{i=1}^{n} \sigma^2(x_{i})|_{0}},
\end{array}
\end{equation}
When the value of $\alpha_{u}^{T}$ at any generation $T$ becomes less than $\alpha_{u}^{stop}$ then the algorithm terminates. In the above equation, $n$ is the number of upper level variables, $\sigma^2(x_{i})|_{T}$ is the variance across dimension $i$ at generation $T$ and $\sigma^2(x_{i})|_{0}$ is the variance across dimension $i$ in the initial population. 
A similar termination scheme can be used when the evolutionary algorithm is executed at the lower level.


\subsection{Lower Level Optimization}\label{sec:lowerLevelAlgorithm}
At the lower level, we utilize SQP if the problem is convex, otherwise we use the lower level evolutionary algorithm described in Table \ref{tab:lowerLevelAlgorithm} that uses similar genetic operations as used at the upper level.

\subsection{Offspring Update}\label{sec:intermediatePhi}
For an offspring $x^{(j)}=(x^{(j)}_u,x^{(j)}_l)$, the lower level vector $x^{(j)}_l$ is updated either using $\Psi$-approximation or $\varphi$-approximation. An update using $\Psi$-approximation is straightforward. However, if an update has to be done using $\varphi$-approximation, it requires to solve the following auxiliary optimization problem. In the auxiliary problem $x_u$ is fixed as $x^{(j)}_u$ and the problem is solved only with respect to $x_l$. The optimal $x_l$ replaces the lower level vector $x^{(j)}_l$ of the offspring.
\begin{align*}
\minimize_{x_l\in X_L}\quad & \hat{F}(x_u,x_l) \\
\st\quad  & \\
 & \hspace{-12mm} \hat{f}(x_u,x_l) \le \hat{\varphi}(x_u) \\
 & \hspace{-12mm} \hat{g}_j(x_u,x_l)\leq 0, j=1,\dots,J\\
 & \hspace{-12mm} \hat{G}_k(x_u,x_l)\leq 0, k=1,\dots,K.
\end{align*}
In the above formulation we use hat for all the functions and constraints as we solve the auxiliary problem on approximated functions and constraints. We use linear approximations for all the constraints, while quadratic approximation is used for the other functions. The auxiliary problem may have to be solved frequently if the lower level problem contains multiple optimal solutions. Solving the auxiliary problem with approximated functions helps in saving actual function evaluations. Note that in the ideal case the auxiliary problem will lead to an optimistic lower level solution corresponding to the fixed $x^{(j)}_u$.

\subsection{Local Search}\label{sec:localSearch}
The algorithm utilizes local search after every $k$ generations of the algorithm. The local search is performed by meta-modeling the upper and lower level functions and constraints along with the $\Psi$ and the $\varphi$-mappings in the vicinity of the best member in the population. Once the $\Psi$ and the $\varphi$-mappings are available, the quality of the two mappings are assessed by the mean square error of the approximations (i.e. $e_{mse}^{\Psi}$ and $e_{mse}^{\varphi}$). The better mapping and the corresponding single level reduction (described in Section \ref{sec:psiMapping} and \ref{sec:varphiMapping}) with approximated functions is solved using SQP to arrive at $x_{u}^{(LS)}$. A lower level optimization corresponding to $x_{u}^{(LS)}$ is solved and if the member is found to be better than $x_{best}^{(j)}$ then $x_{best}^{(j)}$ is updated. In case the member is not better than the best member found so far, then the next local search is performed using the exact upper/lower level objective functions and constraints.

\subsection{Parameters and Platform}\label{sec:architecture}
The algorithm has been implemented in MATLAB. At the upper and lower level, the parameters used in the algorithm are:
\begin{enumerate}
\item $\mu=3$
\item $\lambda=2$
\item $r=2$
\item Probability of crossover = 0.9
\item Probability of mutation = 0.1
\item N = 50 (Population size at upper level)
\item n = 50 (Population size at lower level)
\end{enumerate}

\section{Results}


In this study, we consider three algorithms, the nested approach described in Figure~\ref{fig:nested}, BLEAQ \cite{my-ejor17,my-arxiv13,my-cec14}, and our proposed BLEAQ-II. To assess the performance of each algorithm, 31 runs have been performed for each test instance. During every simulation process, the algorithm is terminated if the function value accuracy reaches the objective function accuracy of $10^{-2}$ at both levels. For each run, the upper and lower level function evaluations required until termination is recorded separately.

\subsection{Standard test problems}

\begin{table*}[htbp]
	\caption{Median function evaluations on TP test suite}
	\begin{center}
		\begin{tabular}{cccc|ccc|cc}
			\toprule
			& \multicolumn{ 3}{c}{UL Func. Evals.} & \multicolumn{ 3}{c}{LL Func. Evals.} & \multicolumn{ 2}{c}{BLEAQ-II Savings} \\ \midrule
			& BLEAQ-II & BLEAQ & Nested & BLEAQ-II & BLEAQ & Nested & $\frac{\mbox{BLEAQ} - \mbox{BLEAQ-II}}{\mbox{BLEAQ}}$ & $\frac{\mbox{Nested} - \mbox{BLEAQ-II}}{\mbox{Nested}}$ \\
			& Med & Med & Med & Med & Med & Med & & \\ \midrule
TP1 & 136 & 155 & - & 242 & 867 & - & 63\% & 98\% \\
TP2 & 255 & 185 & 436 & 440 & 971 & 5,686 & 40\% & 63\% \\
TP3 & 158 & 155 & 633 & 224 & 894 & 6,867 & 64\% & 98\% \\
TP4 & 198 & 357 & 1,755 & 788 & 1,772 & 19,764 & 54\% & 98\% \\
TP5 & 272 & 243 & 576 & 967 & 1,108 & 6,558 & 8\% & 84\% \\
TP6 & 161 & 155 & 144 & 323 & 687 & 1,984 & 43\% & 97\% \\
TP7 & 112 & 255 & 193 & 287 & 987 & 2,870 & 68\% & 95\% \\
TP8 & 241 & 189 & 403 & 467 & 913 & 7,996 & 36\% & 61\% \\ \bottomrule
		\end{tabular}
	\end{center}
	\label{tab:tpResults}
\end{table*}

We first present the empirical results on 8 standard test problems selected from the literature (referred to as TP1-TP8). The description for these test problems has been provided in the Appendix \ref{appendix:TPall}. Table~\ref{tab:tpResults} contains the median upper level (UL) function evaluations, lower level (LL) function evaluations and BLEAQ-II's overall function evaluation savings as compared to other approaches from 31 runs of the algorithms. The overall function evaluations for any algorithm is simply the sum of upper and lower level function evaluations. For instance, for the median run with TP1, BLEAQ-II requires $63\%$ less overall function evaluations as compared to BLEAQ, and $98\%$ less overall function evaluations as compared to the nested approach.

\begin{figure*}
\begin{minipage}[t]{0.49\linewidth}
\begin{center}
\epsfig{file=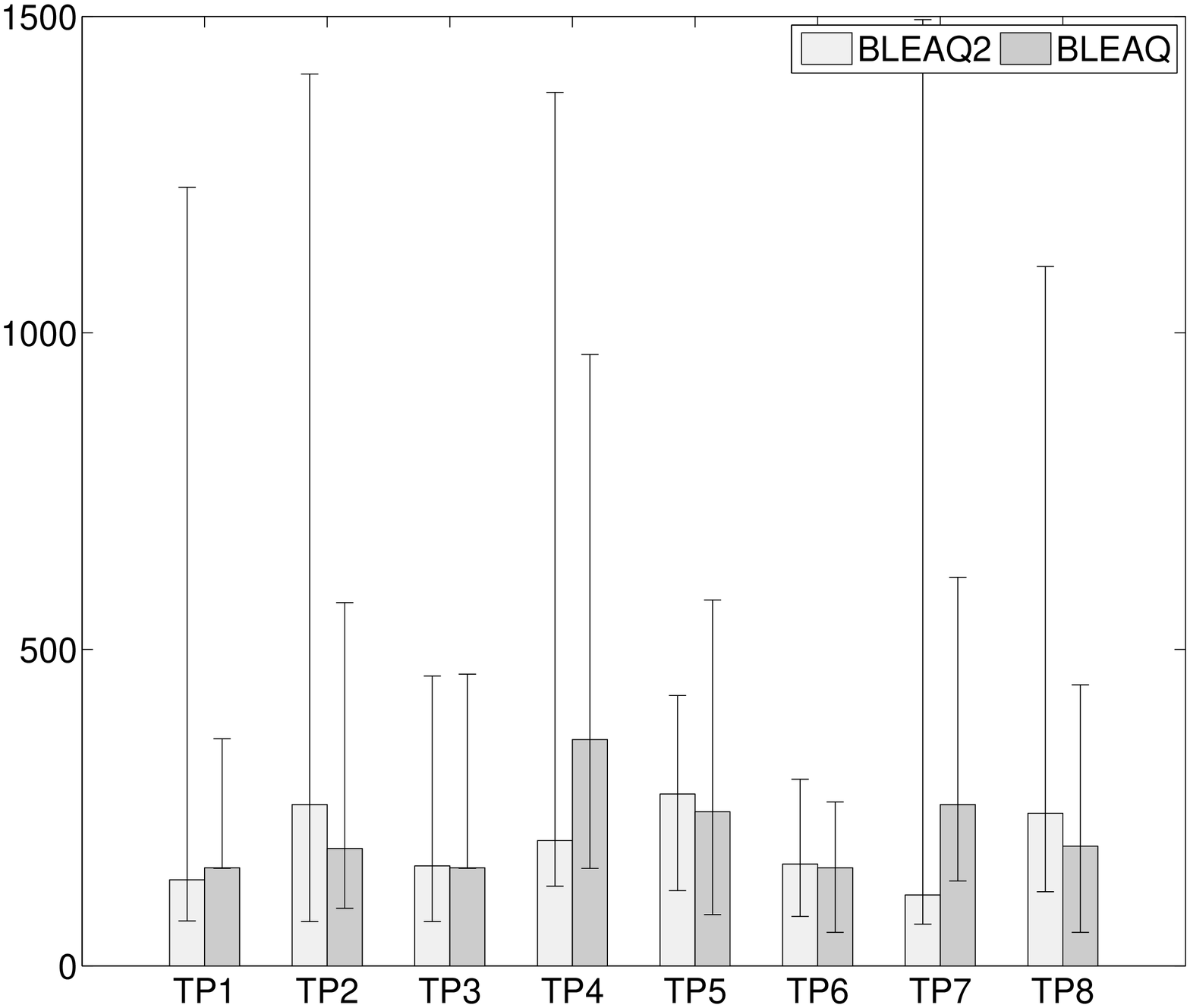,width=0.9\linewidth}
\end{center}
\caption{Bar chart (31 runs/samples) for the upper level function evaluations required for TP 1 to 8.}
\label{fig:tpFEupper}
\end{minipage}\hfill
\begin{minipage}[t]{0.49\linewidth}
\begin{center}
\epsfig{file=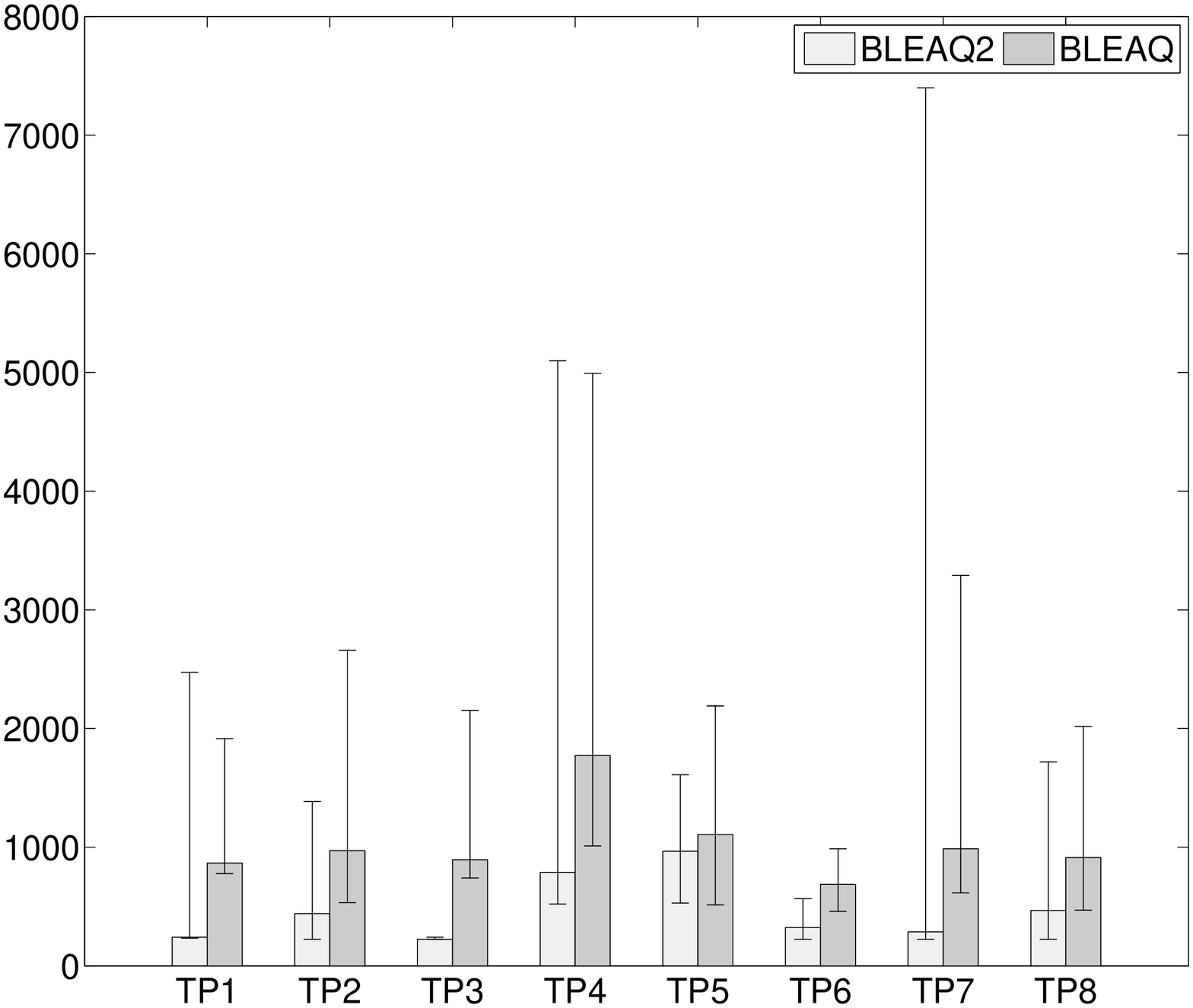,width=0.9\linewidth}
\end{center}
\caption{Bar chart (31 runs/samples) for the lower level function evaluations required for TP 1 to 8.}
\label{fig:tpFElower}
\end{minipage}
\end{figure*}

All these test problems are bilevel problems with small number of variables, and all the three algorithms were able to solve the 8 test instances successfully. A significant computational saving can be observed for both BLEAQ-II and BLEAQ, as compared to the nested approach as shown in the Savings column of Table~\ref{tab:tpResults}. The performance gain going from BLEAQ to BLEAQ-II is quite significant for these simple test problems even though none of them lead to multiple lower level optimal solutions. Detailed comparison between BLEAQ and BLEAQ-II in terms of upper and lower level function evaluations is provided through Figures~\ref{fig:tpFEupper} and \ref{fig:tpFElower}.

\subsection{Scalable test problems}
Next, we compare the results for the three algorithms on the scalable SMD test suite which contains 12 test problems in the original paper \cite{my-ecj14}. We extend this test suite in this paper to a set of 14 test problems by adding two additional scalable test problems. The description for the additional SMD test problems can be found in Appendix \ref{appendix:SMDnew}. First we analyze the performance of the algorithms on a smaller version of the test problems which consists of 5 variables, and then we provide the comparison results on 10-variable instances of the SMD test problems. For the 5 variable version of the SMD test problems, we used the settings as $p=1$, $q=2$ and $r=1$ for all SMD problems except SMD6 and SMD14. For the 5 variable version of SMD6 and SMD14, we used $p=1$, $q=0$, $r=1$ and $s=2$. For the 10 variable version of the SMD test problems, we used the settings as $p=3$, $q=3$ and $r=2$ for all SMD problems except SMD6 and SMD14. For the 10 variable version of SMD6 and SMD14, we used $p=3$, $q=1$, $r=2$ and $s=2$.

\begin{figure*}[hbt]
	\vspace{-4mm}
	\begin{center}
		\epsfig{file=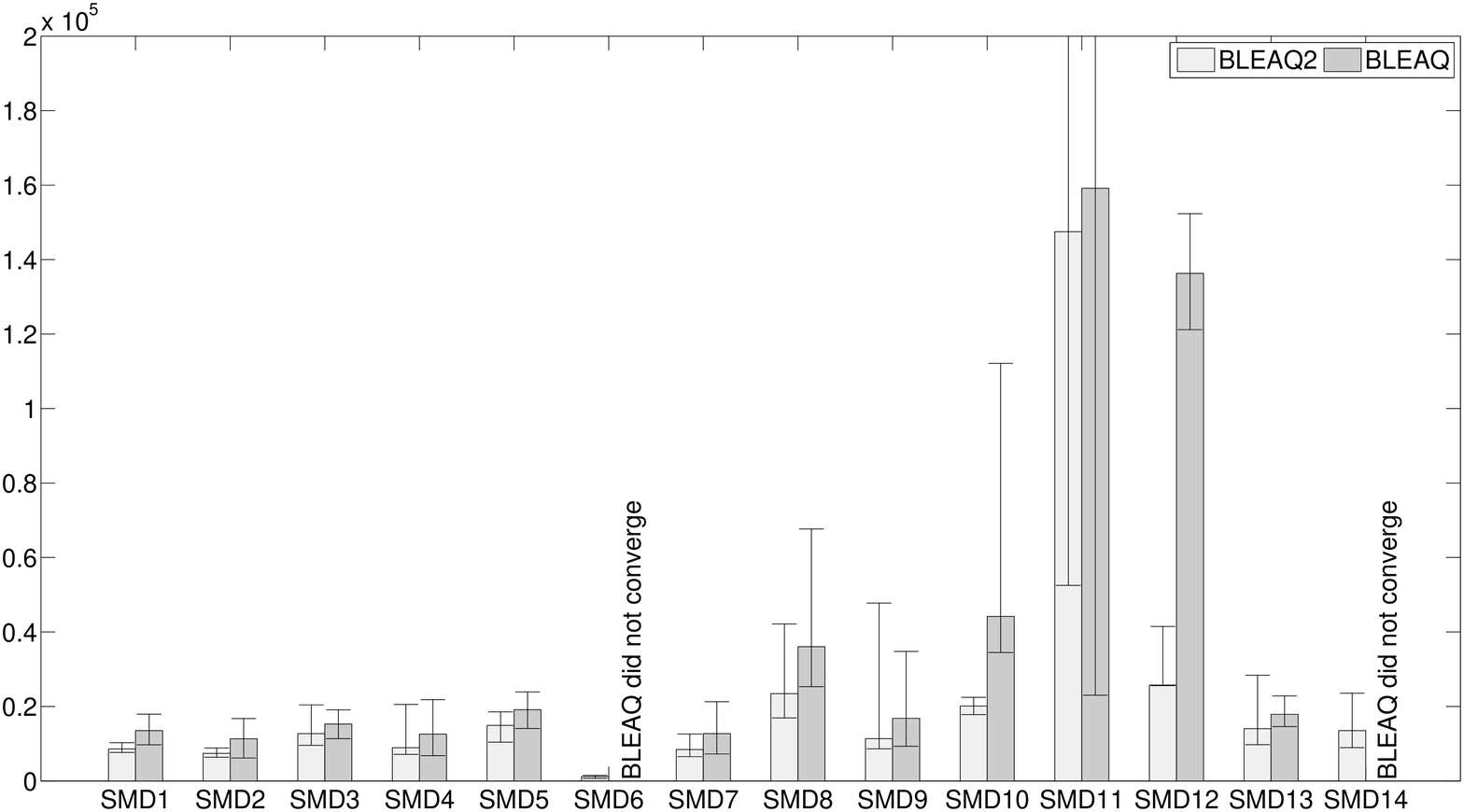,width=.8\linewidth}
	\end{center}
	\vspace{-5mm}
	\caption{Bar chart for overall function evaluations for SMD 1 - 14.}
	\label{fig:smdFE}
	\vspace{-2mm}
\end{figure*}

Table \ref{tab:smdLowResults} provides the median function evaluations and overall savings for the three algorithms on the set of 14 SMD problems. These test problems contain 2 variables at the upper level and 3 variables at the lower level and offer a variety of tunable complexities to the algorithms. For instances, the test set contains problems which are multimodal at the upper and the lower levels, contain multiple optimal solutions at the lower level, contain constraints at the upper and/or lower levels etc. It can be found that BLEAQ-II is able to solve the entire set of 14 SMD test problems, while BLEAQ fails on 2 test problems. The overall savings with BLEAQ-II is higher as compared to BLEAQ for all the test problems. The test problems that contain multiple lower level solutions include SMD6 and SMD14, for which BLEAQ is unable to handle the problem. Further details about the required overall function evaluations from 31 runs can be found in Figures \ref{fig:smdFE}.

\begin{table*}[htbp]
\caption{Median function evaluations on low dimension SMD test suite }
\begin{center}
\begin{tabular}{cccc|ccc|cc}
\toprule
 & \multicolumn{ 3}{c}{UL Func. Evals.} & \multicolumn{ 3}{c}{LL Func. Evals.} & \multicolumn{ 2}{c}{BLEAQ-II Savings} \\ \midrule
 & BLEAQ-II & BLEAQ & Nested & BLEAQ-II & BLEAQ & Nested & $\frac{\mbox{BLEAQ} - \mbox{BLEAQ-II}}{\mbox{BLEAQ}}$ & $\frac{\mbox{Nested} - \mbox{BLEAQ-II}}{\mbox{Nested}}$ \\
 & Med & Med & Med & Med & Med & Med &  &  \\ \midrule
SMD1 & 123  & 98 & 164 & 8,462  & 13,425 & 104,575 & 37\% & 92\% \\
SMD2 & 114  & 88 & 106 & 7,264  & 11,271 & 74,678 & 35\% & 90\% \\
SMD3 & 264  & 91 & 136 & 12,452  & 15,197 & 101,044 & 17\% & 87\% \\
SMD4 & 272  & 110 & 74 & 8,600  & 12,469 & 59,208 & 29\% & 85\% \\
SMD5 & 126  & 80 & 93 & 14,490  & 19,081 & 73,500 & 24\% & 80\% \\
SMD6 & 259  & - & 116 & 914  & - & 3,074 & Large & 63\% \\
SMD7 & 180  & 98 & 67 & 8,242  & 12,580 & 56,056 & 34\% & 85\% \\
SMD8 & 644  & 228 & 274 & 22,866  & 35,835 & 175,686 & 35\% & 87\% \\
SMD9 & 201  & 125 & 127 & 10,964  & 16,672 & 101,382 & 34\% & 89\% \\
SMD10 & 780  & 431 & - & 19,335  & 43,720 & - & 54\% & Large \\
SMD11 & 1735  & 258 & 260 & 134,916  & 158,854 & 148,520 & 14\% & 8\% \\
SMD12 & 203  & 557 & - & 25,388  & 135,737 & - & 81\% & Large \\
SMD13 & 317  & 126 & 211 & 13,729  & 17,752 & 138,089 & 21\% & 90\% \\
SMD14 & 1,014  & - & 168 & 12,364  & - & 91,197 & Large & 85\% \\ \bottomrule
\end{tabular}
\end{center}
\label{tab:smdLowResults}
\end{table*}

\begin{figure}[hbt]
\vspace{-4mm}
\begin{center}
\epsfig{file=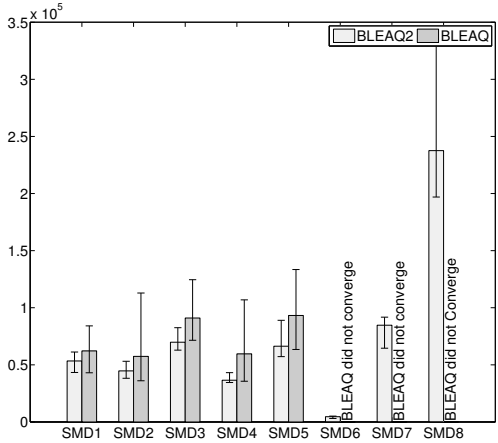,width=0.9\linewidth}
\end{center}
\vspace{-5mm}
\caption{Bar chart for overall function evaluations for 10-dimension SMD 1 - 8.}
\label{fig:smdFEHigh}
\vspace{-2mm}
\end{figure}

Results for the high dimensional SMD test problems have been provided in Table \ref{tab:smdHighResults}. BLEAQ-II leads to much higher savings as compared to BLEAQ, and with higher dimensions BLEAQ is found to once again fail on SMD6 and also on SMD7 and SMD8. Both methods outperform the nested method on most of the test problems. We do not provide results for SMD9 to SMD14 as none of the algorithms were able to handle these problems. It is noteworthy that SMD9 to SMD14 offer difficulties like multi-modalities and highly constrained regions, which none of the algorithms were able to handle with the parameter setting used in this paper. Details for the 31 runs on each of these test problems can be found in Figure \ref{fig:smdFEHigh}.

\begin{figure*}[htp]
\begin{minipage}[t]{0.49\linewidth}
\begin{center}
\epsfig{file=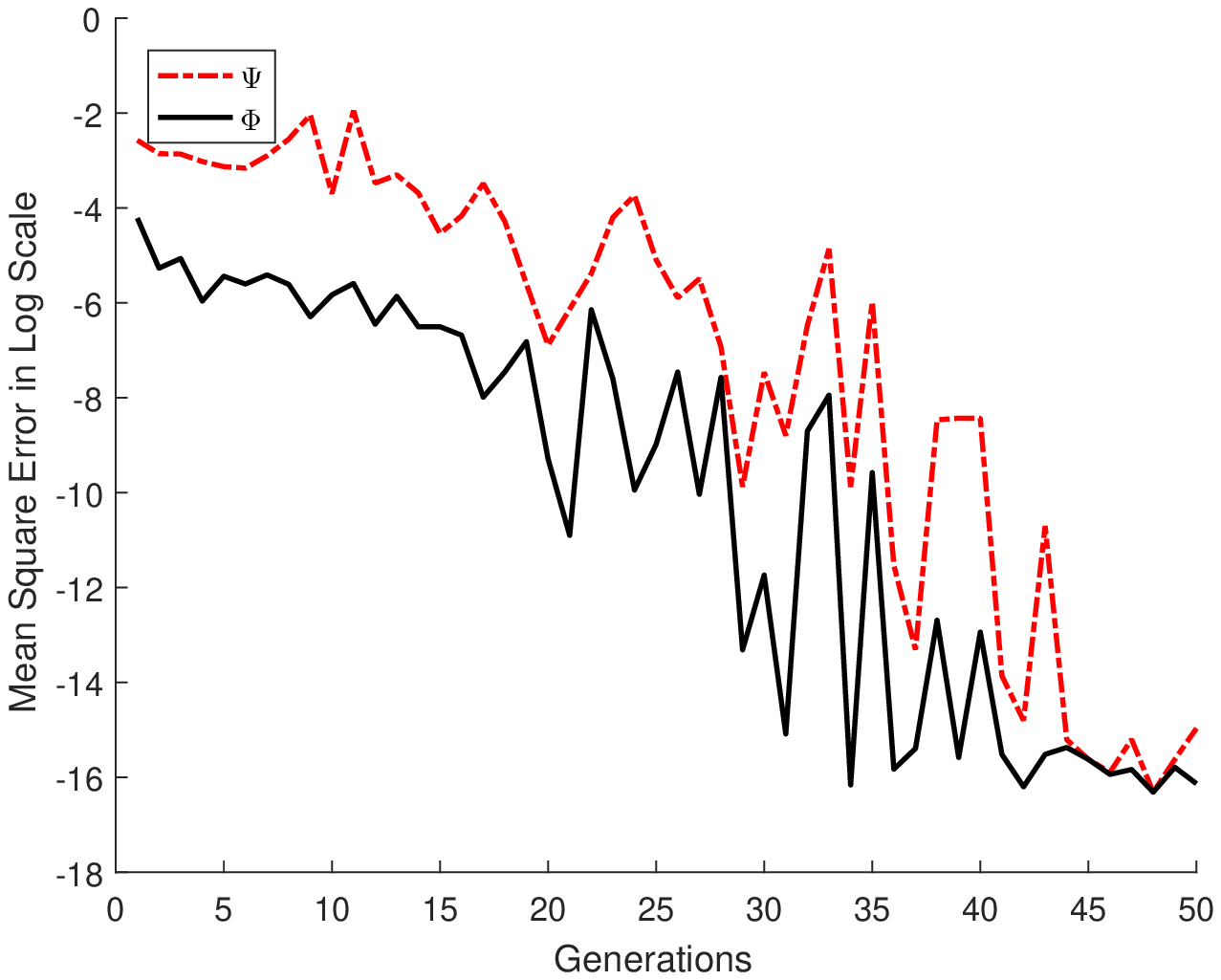,width=0.8\linewidth}
\end{center}
\caption{Approximation error (in terms of Euclidean distance) of a predicted lower level optimal solution when using localized $\Psi$ and $\varphi$-mapping during the intermediate generations of the BLEAQ-II algorithm on the 5-variable SMD1 test problem.}
\label{fig:PsiVsPhiSMD1}
\end{minipage}\hfill
\begin{minipage}[t]{0.49\linewidth}
\begin{center}
\epsfig{file=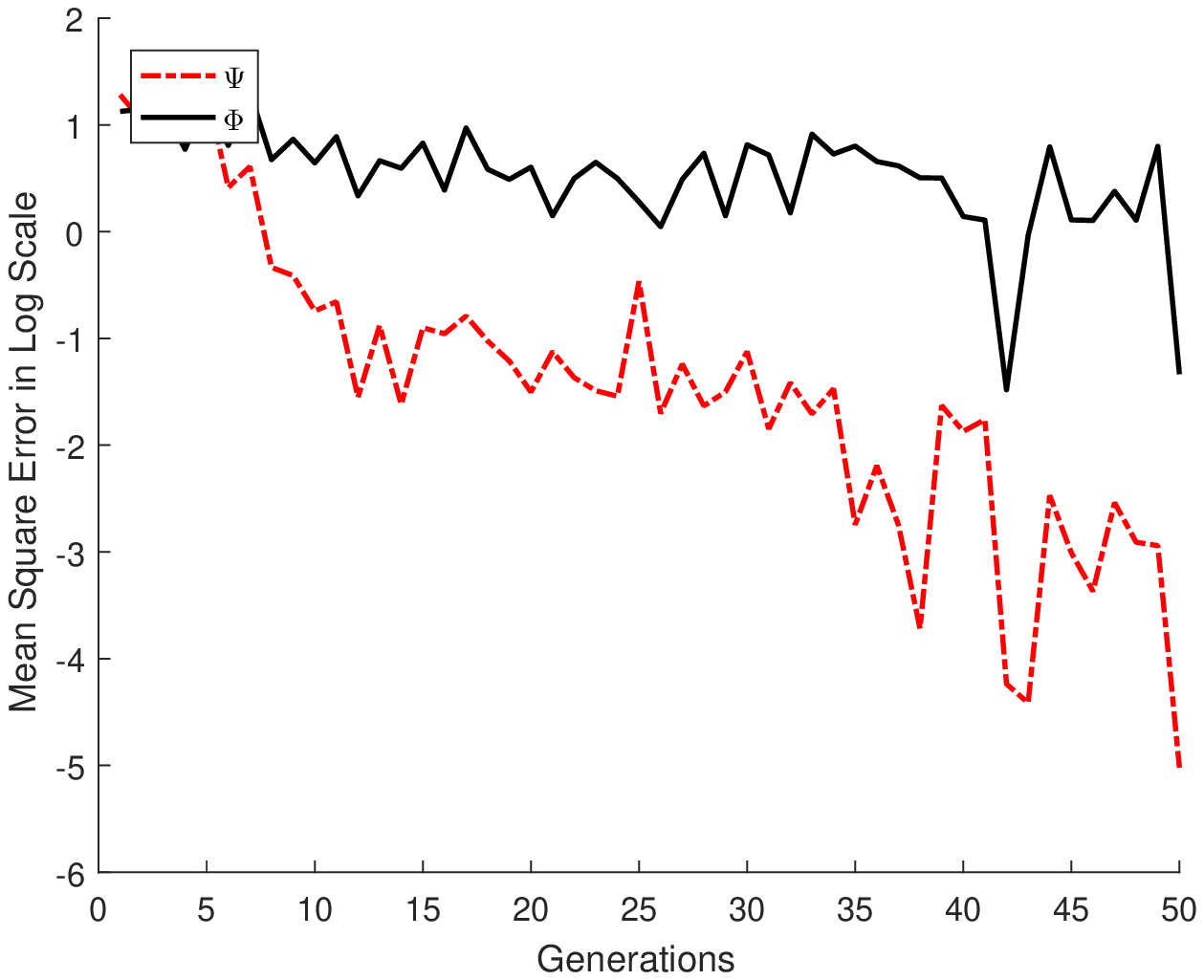,width=0.8\linewidth}
\end{center}
\caption{Approximation error (in terms of Euclidean distance) of a predicted lower level optimal solution when using localized $\Psi$ and $\varphi$-mapping during the intermediate generations of the BLEAQ-II algorithm on the 5-variable SMD13 test problem.}
\label{fig:PsiVsPhiSMD13}
\end{minipage}
\end{figure*}

Through Figures \ref{fig:PsiVsPhiSMD1} and \ref{fig:PsiVsPhiSMD13}, we provide the quality of prediction of the lower level optimal solution made by the $\Psi$-mapping and $\varphi$-mapping approach over the course of the algorithm. It is interesting to note that the quality of $\varphi$-approximation is better in the case of SMD1 test problem in Figure \ref{fig:PsiVsPhiSMD1}, therefore, the prediction decisions are mostly made using the $\varphi$-approximation approach. However, for SMD13 in Figure \ref{fig:PsiVsPhiSMD13}, which involves a difficult $\varphi$-mapping, the prediction decisions are made using the $\Psi$-approximation approach. Both these mappings are found to be improving with an increase in generations of the algorithm. The two figures show the adaptive nature of the BLEAQ-II algorithm in choosing the right approximation strategy based on the difficulties involved in a bilevel optimization problem.

\begin{table*}[htbp]
\caption{Median function evaluations on high dimension SMD test suite }
\begin{center}
\begin{tabular}{cccc|ccc|cc}
\toprule
 & \multicolumn{ 3}{c}{UL Func. Evals.} & \multicolumn{ 3}{c}{LL Func. Evals.} & \multicolumn{ 2}{c}{BLEAQ-II Savings} \\ \midrule
 & BLEAQ-II & BLEAQ & Nested & BLEAQ-II & BLEAQ & Nested & $\frac{\mbox{BLEAQ} - \mbox{BLEAQ-II}}{\mbox{BLEAQ}}$ & $\frac{\mbox{Nested} - \mbox{BLEAQ-II}}{\mbox{Nested}}$ \\
 & Med & Med & Med & Med & Med & Med & & \\ \midrule
SMD1 & 670  & 370 & 760 & 52,866  & 61,732 & 1,776,426 & 14\% & 97\% \\
SMD2 & 510  & 363 & 652 & 44,219  & 57,074 & 1,478,530 & 22\% & 97\% \\
SMD3 & 1369  & 630 & 820 & 68,395  & 90,390 & 1,255,015 & 23\% & 94\% \\
SMD4 & 580  & 461 & 765 & 35,722  & 59,134 & 1,028,802 & 39\% & 96\% \\
SMD5 & 534  & 464 & 645 & 65,873  & 92,716 & 1,841,569 & 29\% & 96\% \\
SMD6 & 584  & - & 824 & 3,950  & - & 156,2003 & Large & 99\% \\
SMD7 & 1,486  & - & - & 83,221  & - & - & Large & Large \\
SMD8 & 6,551  & - & - & 231,040  & - & - & Large & Large \\ \bottomrule
\end{tabular}
\end{center}
\label{tab:smdHighResults}
\end{table*}

\section{Conclusions}
In this paper, we have presented a computationally efficient evolutionary algorithm for solving bilevel optimization problems. The algorithm is based on iterative approximations of two important theoretically motivated mappings; namely, the lower level rational reaction mapping and the lower level optimal value function mapping. The paper discusses about the pros and cons of utilizing these mappings in an evolutionary bilevel optimization algorithm by embedding them in a nested approach. Thereafter, an algorithm is developed that adaptively decides to use one of the mappings during the execution based on the characteristics of the bilevel optimization problem being solved. The proposed algorithm has been tested on a wide variety of bilevel test problems and it has been able to perform significantly better than other approaches in terms of computational requirements.


\onecolumn

\begin{appendices}
\section{Standard Test Problems}\label{appendix:TPall}
In this section, we provide some of the standard bilevel test problems chosen from the literature. Most of these test problems are small with only small number of variables at both levels.

\begin{table*}[!h]
	\caption{Standard test problems TP1-TP5. (Note that $x = x_u$ and $y = x_l$)}
	{\footnotesize
		\begin{tabular}{p{.09\textwidth}  >{$\displaystyle}p{.60\textwidth}<{$}  >{$\displaystyle}p{.20\textwidth}<{$}}
			\toprule
			$\mbox{Problem}$ & \mbox{Formulation} & \mbox{Best Known Sol.} \\
			\toprule
			TP1 \\ $n=2$, $m=2$ &
			\begin{array}{l}
				\underset{{(x,y)}}{\mbox{Minimize}} \hspace{1mm} F(x,y) = (x_{1}-30)^2 + (x_{2}-20)^2 - 20 y_{1} + 20 y_{2}, \\
				\mbox{s.t.}\\
				\quad y \in \underset{(y)}{\argmin}\left\lbrace
				\begin{array}{l}
					f(x,y)=(x_{1}-y_{1})^2+(x_{2}-y_{2})^2\\
					0 \le y_{i} \le 10, \quad i=1,2
				\end{array}
				\right\rbrace, \\
				\quad x_{1}+2x_{2} \ge 30, x_{1}+x_{2} \le 25, x_{2} \le 15\\
			\end{array}
			& \begin{array}{l} \\ \\ \\ \\ \\ \\ \\F = 225.0 \\ f = 100.0
			\end{array}
			\\ \midrule
			TP2 \\ $n=2$, $m=2$ &
			\begin{array}{l}
				\underset{{(x,y)}}{\mbox{Minimize}} \hspace{1mm} F(x,y) = 2 x_{1} + 2 x_{2} - 3 y_{1} - 3 y_{2} - 60, \\
				\mbox{s.t.}\\
				\quad y \in \underset{(y)}{\argmin}\left\lbrace
				\begin{array}{l}
					f(x,y)=(y_{1} - x_{1} + 20)^2+(y_{2} - x_{2} + 20)^2 \\
					x_{1} - 2 y_{1} \ge 10, x_{2} - 2 y_{2} \ge 10\\
					-10 \ge y_{i} \ge 20, \quad i=1,2
				\end{array}
				\right\rbrace, \\
				\quad x_{1}+x_{2}+y_{1}-2 y_{2} \le 40,\\
				\quad 0 \le x_{i} \le 50, \quad i=1,2.
			\end{array}
			& \begin{array}{l} \\ \\ \\ \\ \\ \\ \\F = 0.0 \\ f = 100.0
			\end{array}
			\\ \midrule
			TP3 \\ $n=2$, $m=2$ &
			\begin{array}{l}
				\underset{{(x,y)}}{\mbox{Minimize}} \hspace{1mm} F(x,y) = -(x_{1})^2 - 3 (x_{2})^2 - 4 y_{1} + (y_{2})^2, \\
				\mbox{s.t.}\\
				\quad y \in \underset{(y)}{\argmin}\left\lbrace
				\begin{array}{l}
					f(x,y)=2 (x_{1})^2 + (y_{1})^2 - 5 y_{2} \\
					(x_{1})^2 - 2 x_{1} + (x_{2})^2 - 2 y_{1} + y_{2} \ge -3\\
					x_{2} + 3 y_{1} - 4 y_{2} \ge 4\\
					0 \le y_{i}, \quad i=1,2
				\end{array}
				\right\rbrace, \\
				\quad (x_{1})^2 + 2 x_{2} \le 4,\\
				\quad 0 \le x_{i}, \quad i=1,2
			\end{array}
			& \begin{array}{l} \\ \\ \\ \\ \\ \\ \\F = -18.6787 \\ f = -1.0156
			\end{array}
			\\ \midrule
			TP4 \\ $n=2$, $m=3$ &
			\begin{array}{l}
				\underset{{(x,y)}}{\mbox{Minimize}} \hspace{1mm} F(x,y) = -8 x_{1} - 4 x_{2} + 4 y_{1} - 40 y_{2} - 4 y_{3}, \\
				\mbox{s.t.}\\
				\quad y \in \underset{(y)}{\argmin}\left\lbrace
				\begin{array}{l}
					f(x,y)=x_{1} + 2 x_{2} + y_{1} + y_{2} + 2 y_{3}\\
					y_{2} + y_{3} - y_{1} \le 1\\
					2 x_{1} - y_{1} + 2 y_{2} - 0.5 y_{3} \le 1\\
					2 x_{2} + 2 y_{1} - y_{2} - 0.5 y_{3} \le 1\\
					0 \le y_{i}, \quad i=1,2,3
				\end{array}
				\right\rbrace, \\
				\quad 0 \le x_{i}, \quad i=1,2
			\end{array}
			& \begin{array}{l} \\ \\ \\ \\ \\ \\ \\F = -29.2 \\ f = 3.2
			\end{array}
			\\ \midrule
			TP5 \\ $n=2$, $m=2$ &
						\begin{array}{l}
							\underset{{(x,y)}}{\mbox{Minimize}} \hspace{1mm} F(x,y) = r t(x) x - 3 y_{1} - 4 y_{2} + 0.5 t(y) y, \\
							\mbox{s.t.}\\
							\quad y \in \underset{(y)}{\argmin}\left\lbrace
							\begin{array}{l}
								f(x,y)=0.5 t(y) h y - t(b(x)) y\\
								-0.333 y_{1} + y_{2} - 2 \le 0\\
								y_{1} - 0.333 y_{2} - 2 \le 0\\
								0 \le y_{i}, \quad i=1,2
							\end{array}
							\right\rbrace, \\
							\mbox{where}\\
							\quad h = \left( \begin{array}{cc} 1 & 3\\ 3 & 10\\ \end{array} \right),
							b(x) = \left( \begin{array}{cc} -1 & 2\\ 3 & -3\\ \end{array} \right)x,
							r = 0.1\\
							\quad t(\cdot) \mbox{ denotes transpose of a vector}
						\end{array}
						& \begin{array}{l} \\ \\ \\ \\ \\ \\ \\F = -3.6 \\ f = -2.0
						\end{array}
						\\ \bottomrule
		\end{tabular}
	}
	\vspace{-1mm}
	\label{tab:tpTable1}
\end{table*}
\begin{table*}[!h]
	\caption{Standard test problems TP6-TP8. (Note that $x = x_u$ and $y = x_l$)}
	{\footnotesize
		\begin{tabular}{p{.09\textwidth}  >{$\displaystyle}p{.60\textwidth}<{$}  >{$\displaystyle}p{.20\textwidth}<{$}}
			\toprule
			$\mbox{Problem}$ & \mbox{Formulation} & \mbox{Best Known Sol.} \\
			\toprule
			TP6 \\ $n=1$, $m=2$ &
			\begin{array}{l}
				\underset{{(x,y)}}{\mbox{Minimize}} \hspace{1mm} F(x,y) = (x_{1} - 1)^2 + 2 y_{1} - 2 x_{1}, \\
				\mbox{s.t.}\\
				\quad y \in \underset{(y)}{\argmin}\left\lbrace
				\begin{array}{l}
					f(x,y)=(2 y_{1}-4)^2 +\\ (2 y_{2} - 1)^2 + x_{1} y_{1}\\
					4 x_{1} + 5 y_{1} + 4 y_{2} \le 12\\
					4 y_{2} - 4 x_{1} - 5 y_{1} \le -4\\
					4 x_{1} - 4 y_{1} + 5 y_{2} \le 4\\
					4 y_{1} - 4 x_{1} + 5 y_{2} \le 4\\
					0 \le y_{i}, \quad i=1,2
				\end{array}
				\right\rbrace, \\
				\quad 0 \le x_{1}
			\end{array}
			& \begin{array}{l} \\ \\ \\ \\ \\ \\ \\F = -1.2091 \\ f = 7.6145
			\end{array}
			\\ \midrule
			TP7 \\ $n=2$, $m=2$ &
			\begin{array}{l}
				\underset{{(x,y)}}{\mbox{Minimize}} \hspace{1mm} F(x,y) = -\frac{(x_{1}+y_{1})(x_{2}+y_{2})}{1 + x_{1} y_{1} + x_{2} y_{2}}, \\
				\mbox{s.t.}\\
				\quad y \in \underset{(y)}{\argmin}\left\lbrace
				\begin{array}{l}
					f(x,y)=\frac{(x_{1}+y_{1})(x_{2}+y_{2})}{1 + x_{1} y_{1} + x_{2} y_{2}}\\
					0 \le y_{i} \le x_{i}, \quad i=1,2
				\end{array}
				\right\rbrace, \\
				\quad (x_{1})^2 + (x_{2})^2 \le 100\\
				\quad x_{1} - x_{2} \le 0\\
				\quad 0 \le x_{i}, \quad i=1,2
			\end{array}
			& \begin{array}{l} \\ \\ \\ \\ \\ \\ \\F = -1.96 \\ f = 1.96
			\end{array}
			\\ \midrule
			TP8  \\ $n=2$, $m=2$ &
			\begin{array}{l}
				\underset{{(x,y)}}{\mbox{Minimize}} \hspace{1mm} F(x,y) = |2 x_{1} + 2 x_{2} - 3 y_{1} - 3 y_{2} - 60|, \\
				\mbox{s.t.}\\
				\quad y \in \underset{(y)}{\argmin}\left\lbrace
				\begin{array}{l}
					f(x,y)=(y_{1} - x_{1} + 20)^2 +\\ (y_{2} - x_{2} + 20)^2\\
					2y_{1} - x_{1} + 10 \le 0\\
					2y_{2} - x_{2} + 10 \le 0\\
					-10 \le y_{i} \le 20, \quad i=1,2
				\end{array}
				\right\rbrace, \\
				\quad x_{1} + x_{2} + y_{1} - 2 y_{2} \le 40\\
				\quad 0 \le x_{i} \le 50, \quad i=1,2
			\end{array}
			& \begin{array}{l} \\ \\ \\ \\ \\ \\ \\F = 0.0 \\ f = 100.0
			\end{array}
			\\ \bottomrule
		\end{tabular}
	}
	\vspace{-1mm}
	\label{tab:tpTable2}
\end{table*}

\vspace{10mm}

\section{Additional SMD Test Problems}\label{appendix:SMDnew}
SMD test problems \cite{my-ecj14} are a set of 12 scalable test problems that offer a variety of controllable difficulties to an algorithm. We add two more test problems to the previous test-suite in this paper. Both these problems contain a difficult $\varphi$-mapping, among other difficulties. The upper and lower level functions follow the following structure to induce difficulties due to convergence, interaction, and function dependence between the two levels. The vectors $x_u$ and $x_l$ are further divided into two sub-vectors. The $\varphi$-mapping is defined by the function $f_1$.

\begin{equation}
\begin{array}{l}
F(x_u,x_l) = F_1(x_{u1}) + F_2(x_{l1}) + F_3(x_{u2},x_{l2}) \\
f(x_u,x_l) = f_1(x_{u1}, x_{u2}) + f_2(x_{l1}) + f_3(x_{u2},x_{l2})\\
\mbox{where}\\
     \quad \quad x_u = (x_{u1}, x_{u2}) \quad \mbox{and} \quad x_l = (x_{l1}, x_{l2})
\end{array}
\end{equation}

\begin{table*}[!h]
	\caption{SMD Test Problems. (Note that $(x_{u1}, x_{u2}) = (a,b)$ and $(x_{l1}, x_{l2})=(c,d)$)}
	{\footnotesize
		\begin{tabular}{p{.09\textwidth}  >{$\displaystyle}p{.40\textwidth}<{$}  >{$\displaystyle}p{.40\textwidth}<{$}}
			\toprule
			$\mbox{Problem}$ & \mbox{Formulation} & \mbox{Solution} \\
			\toprule
			SMD13 \\ &
			\begin{array}{l}
F_1 = (a_1 - 1)^2 + \sum_{i=1}^{p-1} \big( (a_i - 1)^2 + (a_{i+1} - (a_{i})^{2})^2 \big),\\
F_2 = -\sum_{i=1}^{q} \sum_{j=1}^{i} (c_{j})^2,\\
F_3 = \sum_{i=1}^{r} \sum_{j=1}^{i} (b_{j})^2 - \sum_{i=1}^{r} (b_i - \log d_{i})^2,\\
f_1 = \sum_{i=1}^{p} \big( |a_i| + 2 |\sin(a_i)|\big),\\
f_2 = \sum_{i=1}^{q} \sum_{j=1}^{i} (c_{j})^2,\\
f_3 = \sum_{i=1}^{r} (b_i - \log d_{i})^2,\\
a_{i} \in [-5,10], \hspace{2mm} \forall \hspace{2mm} i \in \{1,2,\ldots,p\},\\
b_{i} \in [-5,e], \hspace{2mm} \forall \hspace{2mm} i \in \{1,2,\ldots,r\},\\
c_{i} \in [-5,10], \hspace{2mm} \forall \hspace{2mm} i \in \{1,2,\ldots,q\},\\
d_{i} \in (0,10], \hspace{2mm} \forall \hspace{2mm} i \in \{1,2,\ldots,r\}.
			\end{array}
			&
\begin{array}{l}
a_{i} = 1 \; \forall \; i,\\
b_{i} = 0 \; \forall \; i,\\
c_{i} = 0 \; \forall \; i,\\
d_{i} = 1 \; \forall \; i.\\
\end{array}
			\\ \midrule
			SMD14 \\ &
\begin{array}{l}
F_1 = (a_1 - 1)^2 + \sum_{i=1}^{p-1} \big( (a_i - 1)^2 + (a_{i+1} - (a_{i})^{2})^2 \big),\\
F_2 = -\sum_{i=1}^{q} |c_{i}|^{i+1} + \sum_{i=q+1}^{q+s} (c_{i})^{2},\\
F_3 = \sum_{i=1}^{r} i (b_{i})^2 - \sum_{i=1}^{r} |d_{i}|,\\
f_1 = \sum_{i=1}^{p} \left\lfloor a_i \right\rfloor,\\
f_2 = \sum_{i=1}^{q} |c_{i}|^{i+1} + \sum_{i=q+1,i=i+2}^{q+s-1} (c_{i+1} - c_{i})^{2},\\
f_3 = \sum_{i=1}^{r} |(b_{i})^{2} - (d_{i})^{2}|,\\
a_{i} \in [-5,10], \hspace{2mm} \forall \hspace{2mm} i \in \{1,2,\ldots,p\},\\
b_{i} \in [-5,10], \hspace{2mm} \forall \hspace{2mm} i \in \{1,2,\ldots,r\},\\
c_{i} \in [-5,10], \hspace{2mm} \forall \hspace{2mm} i \in \{1,2,\ldots,q+s\},\\
d_{i} \in [-5,10], \hspace{2mm} \forall \hspace{2mm} i \in \{1,2,\ldots,r\}.
\end{array}
			&
\begin{array}{l}
a_{i} = 1 \; \forall \; i,\\
b_{i} = 0 \; \forall \; i,\\
c_{i} = 0 \; \forall \; i,\\
d_{i} = 0 \; \forall \; i.
\end{array}
			\\ \midrule
		\end{tabular}
	}
	\vspace{-1mm}
	\label{tab:smdTable2}
\end{table*}

\end{appendices}

\end{document}